\documentclass[11pt]{article}
\usepackage{graphicx,wrapfig}
\usepackage{flafter}
\usepackage{leftidx}
\usepackage{mathrsfs}
\usepackage{amssymb,amsmath}
\usepackage{bbding}
\usepackage{fancyhdr}
\usepackage{mathrsfs}
\usepackage{color}
\usepackage{stmaryrd}
\usepackage{amsfonts}
\usepackage{latexsym}
\usepackage{psfrag}
\usepackage{graphicx,subfigure}
\usepackage{fancyhdr,graphicx}
\usepackage{multicol}
\usepackage{dsfont}
\usepackage{bbm}
\usepackage{booktabs}
\usepackage[center]{caption2}
\usepackage{cite}
\usepackage{epstopdf}
\usepackage{diagbox}
\usepackage{booktabs}
\usepackage [latin1]{inputenc}
\usepackage{multirow}
\usepackage{enumerate}
\usepackage{enumitem}
\usepackage{algpseudocode,algorithm,algorithmicx}
\usepackage{appendix}

\allowdisplaybreaks
\renewcommand{\leq}{\leqslant}

\date{}
\newtheorem{theorem}{Theorem}[section]
\newtheorem{lemma}{Lemma}[section]
\newtheorem{remark}{Remark}
\newtheorem{example}{Example}[section]

\numberwithin{equation}{section}

\newcommand{\zd}{\,\mathrm{d}}
\newcommand{\diff}{\triangledown_{\tau}}
\newcommand{\myvec}[1]{\boldsymbol{#1}}
\newcommand{\abs}[1]{\left|#1\right|}
\newcommand{\absb}[1]{\big|#1\big|}

\newcommand{\bra}[1]{\left(#1\right)}
\newcommand{\brab}[1]{\big(#1\big)}
\newcommand{\braB}[1]{\Big(#1\Big)}
\newcommand{\brat}[1]{(#1)}
\newcommand{\kbra}[1]{\left[#1\right]}
\newcommand{\kbrab}[1]{\big[#1\big]}

\newcommand{\myinner}[1]{\left\langle#1\right\rangle}
\newcommand{\myinnerb}[1]{\big\langle#1\big\rangle}
\newcommand{\myinnerB}[1]{\Big\langle#1\Big\rangle}
\newcommand{\mynorm}[1]{\left\|#1\right\|}
\newcommand{\mynormb}[1]{\big\|#1\big\|}

\newcommand{\timenorm}[1]{\absb{\!\absb{\!\absb{#1}\!}\!}}

\begin{document}
\title{An adaptive BDF2 implicit time-stepping method \\for the phase field crystal model}
\author{
Hong-lin Liao\thanks{ORCID 0000-0003-0777-6832; Department of Mathematics,
Nanjing University of Aeronautics and Astronautics,
Nanjing 211106, P. R. China. Hong-lin Liao (liaohl@csrc.ac.cn, liaohl@nuaa.edu.cn)
is supported by a grant 1008-56SYAH18037
from NUAA Scientific Research Starting Fund of Introduced Talent.}
\quad Bingquan Ji\thanks{Department of Mathematics, Nanjing University of Aeronautics and Astronautics,
211101, P. R. China. Bingquan Ji (jibingquanm@163.com).}
\quad Luming Zhang\thanks{Department of Mathematics, Nanjing University of Aeronautics and Astronautics,
211101, P. R. China. Luming Zhang (zhanglm@nuaa.edu.cn)
is supported by the NSFC grant No. 11571181.}}
\date{}
\maketitle
\normalsize

\begin{abstract}
  An  adaptive BDF2 implicit time-stepping method is analyzed for the  phase field crystal model.
  The suggested method is proved to preserve a modified energy dissipation law at the discrete levels
  if the time-step ratios $r_k:=\tau_k/\tau_{k-1}<3.561$,
  a recent zero-stability restriction of variable-step BDF2 scheme for ordinary differential problems.
  By using the discrete orthogonal convolution kernels and the corresponding convolution inequalities,
  an optimal $L^2$ norm error estimate is established
  under the weak step-ratio restriction $0<r_k<3.561$ ensuring the energy stability.
  This is the first time such error estimate is theoretically proved for a nonlinear parabolic equation.
  On the basis of ample tests on random time meshes,
  a useful adaptive time-stepping strategy is suggested to efficiently capture
  the multi-scale behaviors and to accelerate the numerical simulations.
  \noindent{\bf AMS subject classiffications.}\;\; 35Q99, 65M06, 65M12, 74A50
\end{abstract}

\section{Introduction}\setcounter{equation}{0}

The phase field crystal (PFC) growth model \cite{elder2002modeling}
is an efficient approach to
simulate crystal dynamics at the atomic scale in space while on diffusive scales in time.
This model has been  successfully applied to
a wide variety of simulations in the microstructure evolution \cite{elder2002modeling},
epitaxial thin film growth \cite{elder2004modeling} and materials science across different
time scales \cite{provatas2007using,asadi2015a}.
The phase variable of PFC model describes a coarse-grained temporal average of the number
density of atoms, and the model is thermodynamically consistent in that the free energy of the thermodynamic model is dissipative.
Consider a free energy functional of Swift-Hohenberg type \cite{elder2002modeling,elder2004modeling},
\begin{align}\label{cont:free energy}
E[\Phi] = \int_{\Omega}\bra{\frac{1}{4}\Phi^4
+\frac{1}{2}\Phi\bra{-\epsilon + (1+\Delta)^2}\Phi}\zd\mathbf{x},
\end{align}
where $\mathbf{x}\in\Omega\subseteq\mathbb{R}^d$ ($d=1,2,3$),  $\Phi$ represents the
atomistic density field and $\epsilon\in(0,1)$ is a parameter related to the temperature.
Then the phase field crystal equation is given by the $H^{-1}$ gradient flow
associated with the free energy functional $E[\phi]$,
\begin{align}\label{cont: Problem-PFC}
\partial_t \Phi=\Delta\mu\quad\text{with}\quad
\mu=\tfrac{\delta E}{\delta \Phi}= \Phi^3-\epsilon\Phi+(1+\Delta)^2\Phi,
\end{align}
where $\mu$ is called the chemical potential.
We assume that $\Phi$ is periodic over the domain $\Omega$.
By applying the integration by parts, one can find the volume conservation, $\brab{\Phi(t),1}=\brab{\Phi(t_0),1}$,
and the following energy dissipation law,
\begin{align}\label{cont:energy dissipation}
\frac{\zd{E}}{\zd{t}}=\brab{\tfrac{\delta E}{\delta \Phi},\partial_t \Phi}
=\bra{\mu,\Delta\mu}
=-\mynorm{\nabla\mu}^2 \le 0,
\end{align}
where the $L^2$ inner product $\bra{f,g}:=\int_{\Omega}fg\zd{\mathbf{x}}$,
and the associated $L^2$ norm $\mynorm{f}:=\sqrt{\bra{f,f}}$ for
all $f,g\in L^{2}(\Omega)$.

The PFC equation is a sixth-order nonlinear partial differential equation
and it may be challenging to design efficient and stable numerical algorithms.
As for the time integration approaches,
Crank-Nicolson (CN) schemes
\cite{baskaran2013convergence,Dong2018Convergence,li2019efficient,
Liu2019Efficient,jing2019linear,zhang2013an,yang2017linearly}
and backward differentiation formulas (BDF)
\cite{wise2009an,wang2011an,li2019efficient,
Liu2019Efficient,li2017an,yang2017linearly,Yan2018A,Cheng2019Fourth,Cheng2019An}
are wide-spread in the literatures.
Due to the energy dissipation property \eqref{cont:energy dissipation},
BDF1 and BDF2 methods seem to be more suitable than CN type schemes
in resolving this stiff problem.
Actually, the BDF1 and BDF2 methods both are  A-stable and L-stable,
while the trapezoidal formula is only A-stable.
Moreover, the preservation of \eqref{cont:energy dissipation} at the discrete time levels,
called energy stability, has been regarded as a basic requirement of numerical methods
to be effective in simulating the long-time coarsening dynamics.

The main goal of the existing techniques is to guarantee the energy stability,
including linearized treatments
\cite{Xu2006Stability,Shen2010Numerical,Shen2016On,
Gong2019Energy,li2019efficient,Liu2019Efficient}
and the nonlinear progressing
\cite{wise2009an,wang2011an,baskaran2013convergence,zhang2013an,li2017an}.
The linearized treatments always lead to a linear system of algebraic equations,
which improve the computational efficiency since they avoid an inner iteration.
There are many linearized strategies, such as
the stabilized methods \cite{Xu2006Stability,Shen2010Numerical,Shen2016On},
the invariant energy quadratization (IEQ) method
\cite{Gong2019Energy,li2019efficient,Liu2019Efficient},
and the scalar auxiliary variable (SAV) approach
\cite{Gong2019Energy,li2019efficient,Liu2019Efficient}.
Precisely, the stabilized semi-implicit methods use some appropriate high-order
linear terms to construct linearly energy stable schemes.
The common goal of IEQ and SAV methods is to transform the original system into
a new equivalent system with a quadratic energy functional preserving
the corresponding modified energy dissipation property.
We note that SAV approach usually leads to numerical schemes
involving only the decoupled equations with constant coefficients.
As is known to all, the linearized treatments require small time steps to
control the linearization error or ensure the stability.
However, large time steps are necessary to accelerate the numerical simulations,
especially in the coarsening process of phase field models.

In recent years, the nonlinear treatments,
mainly involving the convex splitting techniques
\cite{wise2009an,Dong2018Convergence,wang2011an,baskaran2013convergence}
and fully implicit methods \cite{zhang2013an,li2017an},
have also received extensive attentions.
In the framework of convex splitting strategy,
the convex and concave parts of chemical potential
are treated implicity and explicitly, respectively.
It results in a nonlinear scheme having the unique solvability
and unconditionally energy stability.
As pointed out by Xu et al. \cite{xu2019on},
a major advantage of convex splitting implicit schemes is that a relatively
large time-step size can be used;
but such schemes with large time-step sizes may have time delays
and hence may be inaccurate.
Actually, the convex splitting scheme
can be mathematically interpreted as a full implicit scheme of a convexified model with a
time-delay regularized term of the original equation,
see more details in \cite{xu2019on}.
Numerical evidences indicate that the convex splitting techniques usually
lead to approximation of the solution of the original model at a delayed time,
especially when large time-steps are used.
So the fully implicit schemes are recommended by Xu et al. \cite{xu2019on}
since they are workable for large time-step sizes,
and avoid the potential time-delays in the long-time numerical simulations.

As a remarkable feature of phase field problems including the PFC equation,
they always permit multiple time scales in approaching the steady state.
Therefore, the adaptive time-stepping strategy would be much
more preferred to  resolve varying time scales efficiently
and to reduce the computational cost significantly.
In the literature, some commonly used adaptive time-stepping strategies
consist of utilizing
the accuracy criterion \cite{Gomez2011Provably}
and the time derivative  of the total energy \cite{Qiao2011An,zhang2013an}.
More precisely, the adaptive time step method reported in \cite{Gomez2011Provably} permits
large time steps when the solution is smooth,
and uses small time steps when the solution is less regular.
The adaptive technique used in \cite{Qiao2011An} produces small time steps
when the energy decays rapidly,
and permits large time steps when the energy decays slowly.
Considerable numerical evidences
showed that both of them can greatly save the computational cost.

This paper considers an adaptive BDF2 implicit time-stepping method
for the PFC equation.  Consider the nonuniform time levels
$0=t_{0}<t_{1}<\cdots<t_N=T$ with the time-step sizes
$\tau_{k}:=t_{k}-t_{k-1}$ for $1\le k \le N$,
and denote the maximum time-step size $\tau:=\max_{1\le k\le N}\tau_k$.
Let the local time-step ratio $r_k:=\tau_k/\tau_{k-1}$ for $2\le k\le N$,
and let $r_1\equiv0$ when it appears.
Given a grid function $\{v^k\}_{k=0}^N$,
put $\diff v^{k}:=v^{k}-v^{k-1}$, $\partial_{\tau}v^{k}:=\diff v^{k}/\tau_k$ for $k\geq{1}$.
Taking $v^n=v(t_n)$, we always view the variable-step BDF2 formula as a discrete convolution
summation
\begin{align}\label{def: BDF2-Formula}
D_2v^n:=\sum_{k=1}^nb_{n-k}^{(n)}\diff v^k
\quad \text{for $n\ge1$},
\end{align}
in which the discrete convolution kernels $b_{n-k}^{(n)}$ are
defined by, for $n\ge 2$,
\begin{align}\label{def:BDF2-kernels}
b_{0}^{(n)}:=\frac{1+2r_n}{\tau_n(1+r_n)},\quad
b_{1}^{(n)}:=-\frac{r_n^2}{\tau_n(1+r_n)}\quad \text{and} \quad
b_{j}^{(n)}:=0,\quad \mathrm{for}\quad 2\le j\le n-1,
\end{align}
and $b_0^{(1)}:=1/\tau_1$ when $n=1$.
Obviously, by taking $r_1=0$, the BDF2 scheme \eqref{def: BDF2-Formula} reduces to the BDF1 method for $n=1$.
Here we will use the BDF1 scheme to compute the first-level
numerical solution having the second-order temporal accuracy.

It is known that the rigorous numerical analysis of nonuniform one-step approaches
might be relatively easy since they contain only one degree of freedom,
i.e., the current time step size.
By contrast, the numerical analysis of multi-step methods involving multiple degrees of freedom
(the current and previous time step sizes)
seems rather difficult, especially on a general class of time meshes.
For the underlaying variable-step BDF2 method for ordinary initial-value problems,
Grigorieff \cite{grigorieff1983stability} proved almost forty years ago
that it is zero-stable only if the adjacent time-step ratios $r_k<1+\sqrt{2}$.
Twenty years ago,  Becker \cite{Becker1998} applied the variable-step
BDF2 formula to a linear parabolic equation
and established a second-order temporal convergence
only if $r_k\le(2+\sqrt{13})/3\approx 1.868$.
However, the resulting error estimate is far from sharp
because it involves an undesired prefactor $\exp(C\Gamma_n)$
where $\Gamma_n$ may be unbounded as the time step sizes vanish.
Recently, Chen et al. \cite{chen2019a} analyzed a variable-step
stabilized BDF2 scheme for the Cahn-Hilliard equation.
This work replaced the undesirable prefactor $\exp(C\Gamma_n)$ by
a bounded exponential prefactor $\exp(Ct_n)$
with the help of a generalized Gr\"{o}nwall inequality.
Nonetheless, it seems that the somewhat rigid restriction $r_k\le1.53$ in \cite{chen2019a}
may be hard to weaken due to the combined technique
using the $H^1$ norm error to control the $L^2$ norm error.

Recently, the variable-step BDF2 method was revisited
in our previous report \cite{Liao2019Adaptive} from
a new point of view by making use of the positive
semi-definiteness of BDF2 kernels $b_{n-k}^{(n)}$.
As a result, a concise $L^2$ norm convergence theory of adaptive BDF2
scheme for linear diffusion equation was established
provided  the adjacent time-step ratios $r_k\le(3+\sqrt{17})/2\approx 3.561$.
The main discrete tool used in \cite{Liao2019Adaptive} is the
discrete orthogonal convolution (DOC) kernels, that is,
\begin{align}\label{DOC-Kernels}
\theta_{0}^{(n)}:=\frac{1}{b_{0}^{(n)}}
\quad \mathrm{and} \quad
\theta_{n-k}^{(n)}:=-\frac{1}{b_{0}^{(k)}}
\sum_{j=k+1}^n\theta_{n-j}^{(n)}b_{j-k}^{(j)}
\quad \text{for $1\le k\le n-1$}.
\end{align}
One has the following  discrete orthogonal identity
   \begin{align}\label{orthogonal identity}
 \sum_{j=k}^n\theta_{n-j}^{(n)}b_{j-k}^{(j)}\equiv \delta_{nk}\quad\text{for $1\le k\le n$,}
   \end{align}
where $\delta_{nk}$ is the Kronecker delta symbol. By exchanging the summation order
and using the identity \eqref{orthogonal identity}, it is not difficult to check that
  \begin{align}\label{orthogonal equality for BDF2}
  \sum_{j=1}^n\theta_{n-j}^{(n)}D_2v^j=\diff v^n\quad\text{for any sequence $\{v^j\,|\,0\le j\le n\}$.}
   \end{align}
The equality \eqref{orthogonal equality for BDF2} will play an important role in the subsequent analysis.
More properties of the DOC kernels $\theta_{n-k}^{(n)}$ are referred to Lemma \ref{lem:DOC property} below.

In this paper, we continue to develop the recent technique in \cite{Liao2019Adaptive}
and derive some novel discrete convolution inequalities with respect to the DOC kernels $\theta_{n-k}^{(n)}$.
An optimal $L^2$ error estimate of
the fully implicit BDF2 scheme with unequal time-step sizes
is achieved for solving the PFC equation \eqref{cont: Problem-PFC},
\begin{align}\label{scheme:PFC BDF2}
D_2\phi^n=\Delta_h\mu^n\quad\text{with}\quad
\mu^n=(1+\Delta_h)^2\phi^n + \brat{\phi^n}^3-\epsilon\phi^n
\quad\text{for $1\le n\le N$,}
\end{align}
subject to the periodic boundary conditions  and a proper initial data $\phi^0\approx\Phi^0$.
The spatial operators are approximated by
the Fourier pseudo-spectral method, as described in the next section.
Firstly, the unique solvability is established in Theorem \ref{thm: convexity solvability}
by using the fact that the solution of nonlinear scheme \eqref{scheme:PFC BDF2}
is equivalent to the minimization of a convex functional.
Lemma \ref{lem:conv kernels positive} shows that the BDF2 convolution kernels $b_{n-k}^{(n)}$ are positive definite
provided the adjacent time-step ratios $r_k$ satisfy a sufficient condition
\begin{enumerate}[itemindent=1em]
\item[\textbf{S1}.]
 $0< r_k < r_{\mathrm{sup}}:=\bra{3+\sqrt{17}}/2 \approx 3.561$ for $2\le k \le N$.
\end{enumerate}
We then verify in Theorem \ref{thm:energy decay law} that
the adaptive BDF2 time-stepping method \eqref{scheme:PFC BDF2}
preserves a modified energy dissipation law at the discrete time levels
under a proper step-size restriction.
The maximum norm bound of solution is obtained in Lemma \ref{lem:Bound-Solution}
so that the subsequent error estimate can be derived without
assuming the Lipschitz continuity of nonlinear bulk force.

Section 3 focuses on the $L^2$ norm convergence of the suggested adaptive BDF2 method \eqref{scheme:PFC BDF2}.
The main tools are the above DOC kernels $\theta_{n-k}^{(n)}$ defined in \eqref{DOC-Kernels}
and the corresponding discrete convolution inequalities, see Lemmas \ref{lem:quadr form inequ} and \ref{lem:inner quad ineq}.
Although the condition \textbf{S1} permits us to use a series of increasing time-steps
with the amplification factors up to 3.561,
very large time-steps always result in a loss of numerical
accuracy. So large amplification factors would be rarely appeared continuously
in practice and it is reasonable to assume that
\begin{enumerate}[itemindent=1em]
\item[\textbf{S2}.]
  The time-step ratios $r_k$ are contained in  \textbf{S1},
  but almost all of them less than $1+\sqrt{2}$, or $\abs{\mathfrak{R}}=N_0\ll N$,
  where $\mathfrak{R}$ is an index set  $ \mathfrak{R}:=\{k\,|\,1+\sqrt{2}\le r_k < (3+\sqrt{17})/2\}. $
\end{enumerate}
Potential users would be recommended to take $r_k\in(0,1+\sqrt{2})$
with $N_0=0$ in practical numerical simulations.
Also,
as shown in Theorem \ref{thm:Convergence-Results}
and Remark \ref{remark: S1-S2 accuracy order},
this restriction \textbf{S2} ensures the second-order convergence in time.
Several numerical examples are presented in Section 4 to validate the
accuracy and effectiveness of our method \eqref{scheme:PFC BDF2}.

In summary, our contributions in this paper are three folds:
\begin{enumerate}
  \item An energy dissipation law at the discrete time level
  with a modified energy form is established for the BDF2 implicit method \eqref{scheme:PFC BDF2}
  if the adjacent step ratios $r_k$ satisfy \textbf{S1}.
  It leads to the stability in the maximum norm.
  \item The BDF2 implicit method \eqref{scheme:PFC BDF2} is shown
  to be convergent in the $L^2$ norm under the condition \textbf{S1},
  and the second-order accuracy is achieved if \textbf{S2} holds.
  To the best of our knowledge, this is the first time such
  an optimal $L^2$ norm error estimate of variable-step BDF2 method
  is proved for a nonlinear sixth-order parabolic problem.
 \item Extensive numerical experiments and comparisons to the Crank-Nicolson scheme
  are performed to show the effectiveness of BDF2 time-stepping approach,
  especially when coupled with an adaptive time-stepping strategy.
\end{enumerate}
Throughout this paper, any subscripted $C$, such as $C_u$ and $C_\phi$, denotes a generic positive constant, not necessarily
 the same at different occurrences; while,
 any subscripted $c$, such as $c_\Omega,c_0,c_1$ and $c_2$,
denotes a fixed constant.
Always, the appeared constants are  dependent on the given data
and the solution but independent of the time steps and spatial lengths.

\section{Energy dissipation law and solvability}
\setcounter{equation}{0}

\subsection{Spatial discretization and preliminary results}
For simplicity of presentation, set the spatial domain $\Omega=(0,L)^3$
and consider the uniform length $h_x=h_y=h_z=h:=L/M$ in three spatial directions
for an even positive integer $M$.
We define the discrete grid
 $\Omega_{h}:=\big\{\mathbf{x}_{h}=(ih,jh,kh)\,|\,1\le i,j,k \le M\big\}$
and put
$\bar{\Omega}_{h}:=\Omega_{h}\cup\partial\Omega$.
Denote the space of $L$-periodic grid functions
$\mathbb{V}_{h}:=\{v\,|\,v=\bra{v_h}\; \text{is $L$-periodic for}\; \mathbf{x}_h\in\bar{\Omega}_h\}.$
For any grid functions $v,w\in\mathbb{V}_{h}$,
define the discrete inner product
$\myinner{v,w}:=h^3\sum_{\mathbf{x}_h\in\Omega_{h}}v_h w_h$,
the associated $L^{2}$ norm $\mynorm{v}:=\sqrt{\myinner{v,v}}$.
Also, we will use the discrete $L^4$ norm $\mynorm{v}_{l^4}=\sqrt[4]{h^3\sum_{\mathbf{x}_h\in\Omega_{h}}|v_h|^4}$ and the maximum norm $\mynorm{v}_{\infty}:=\max_{\mathbf{x}_{h}\in\Omega_{h}}|v_{h}|$.

For a periodic function $v(\mathbf{x})$ on $\bar{\Omega}$,
let $P_M:L^2(\Omega)\rightarrow \mathscr{F}_M$
be the standard $L^2$ projection operator onto the space $\mathscr{F}_M$,
consisting of all trigonometric polynomials of degree up to $M/2$,
and $I_M:L^2(\Omega)\rightarrow \mathscr{F}_M$
be the trigonometric interpolation operator \cite{Shen2011Spectral},
that is,
\[
\bra{P_Mv}(\mathbf{x})=\sum_{\ell,m,n =- M/2}^{M/2-1}
\widehat{v}_{\ell,m,n}e_{\ell,m,n}(\mathbf{x}),\quad
\bra{I_Mv}(\mathbf{x})=\sum_{\ell,m,n =- M/2}^{M/2-1}
\widetilde{v}_{\ell,m,n}e_{\ell,m,n}(\mathbf{x}),
\]
where the complex exponential basis functions
$e_{\ell,m,n}(\mathbf{x}):=e^{\mathrm{i}\nu\bra{\ell x+my+nz}}$ with $\nu=2\pi/L$.
The coefficients $\widehat{v}_{\ell,m,n}$
refer to the standard Fourier coefficients of function $v(\mathbf{x})$,
and the
pseudo-spectral coefficients $\widetilde{v}_{\ell,m,n}$ are determined such that $\bra{I_Mv}(\mathbf{x}_h)=v_h$.

The Fourier pseudo-spectral first and second order derivatives of $v_h$ are given by
\[
\mathcal{D}_xv_h:=\sum_{\ell,m,n= -M/2}^{M/2-1}
\bra{\mathrm{i}\nu\ell}\widetilde{v}_{\ell,m,n}
e_{\ell,m,n}(\mathbf{x}_h),\quad
\mathcal{D}_x^2v_h:=\sum_{\ell,m,n = -M/2}^{M/2-1}
\bra{\mathrm{i}\nu\ell}^2\widetilde{v}_{\ell,m,n}
e_{\ell,m,n}(\mathbf{x}_h).
\]
The differentiation  operators $\mathcal{D}_y,\mathcal{D}_y^2,\mathcal{D}_z$
and $\mathcal{D}_z^2$ can be defined in the similar fashion.
In turn, we can define the discrete gradient $\nabla_h$
and Laplacian $\Delta_h$
in the point-wise sense, by
\[
\nabla_hv_h := \left(\mathcal{D}_xv_h,\mathcal{D}_yv_h,\mathcal{D}_zv_h\right)^T\quad\text{and}\quad
\Delta_hv_h :=\nabla_h\cdot\bra{\nabla_hv_h}= \bra{\mathcal{D}_x^2+\mathcal{D}_y^2+\mathcal{D}_z^2}v_h.
\]
For any periodic grid functions $v,w\in\mathbb{V}_{h}$, it is easy to check the following discrete Green's formulas,
see \cite{gottlieb2012stability,cheng2016a} for more details, $\myinner{-\Delta_hv,w}=\myinner{\nabla_hv,\nabla_hw}$,
$\myinner{\Delta_h^2v,w}=\myinner{\Delta_hv,\Delta_hw}$,  and $\myinner{\Delta_h^3v,w}=-\myinner{\nabla_h\Delta_hv,\nabla_h\Delta_hw}$.
Also we have the following embedding inequality
\begin{align}\label{ieq: Linfty H2 embedding}
\mynormb{v}_\infty\le c_\Omega\bra{\mynormb{v}+\mynormb{\Delta_hv}}\quad \text{for any $v\in\mathbb{V}_{h}$.}
\end{align}

For the underlying volume-conservative problem,
it is convenient to define a mean-zero space
$$\mathbb{\mathring V}_{h}:=\big\{v\in\mathbb{V}_{h}\,|\, \myinner{v,1}=0\big\}\subset\mathbb{V}_{h}.$$
As usual, one can introduce a discrete version of
inverse Laplacian operator $\bra{-\Delta_h}^{-\gamma}$
by following the arguments in \cite{cheng2016a}.
For a grid function $v\in\mathbb{\mathring V}_{h}$, define
\[
\bra{-\Delta_h}^{-\gamma}v_h
:=\sum_{\mbox{\tiny$\begin{array}{c}
\ell,m,n=-M/2\\
\bra{\ell,m,n}\neq \mathbf{0}
\end{array}$}}^{M/2-1}
\bra{\nu^2\bra{\ell^2+m^2+n^2}}^{-\gamma}\widetilde{v}_{\ell,m,n}
e_{\ell,m,n}(\mathbf{x}_h),
\]
and an $H^{-1}$ inner product
\[
\myinner{v,w}_{-1}
:=\myinnerb{\bra{-\Delta_h}^{-1}v,w}.
\]
The associated $H^{-1}$ norm $\mynorm{\cdot}_{-1}$ can be defined by $\mynorm{v}_{-1}:=\sqrt{\myinner{v,v}_{-1}}\,.$
We have the following generalized H\"{o}lder inequality,
\begin{align}\label{ieq:generalized Holder}
\mynormb{v}^2\le \mynormb{\nabla_hv}\mynormb{v}_{-1}\quad \text{for any $v\in\mathbb{\mathring V}_{h}$.}
\end{align}

\subsection{Unique solvability}

\begin{lemma}\label{lem: L2-Embedding-Inequality}
For any $v\in \mathbb{\mathring V}_{h}$, it holds that
$\mynormb{v}^2 \le \frac{1}{3}\mynormb{(1+\Delta_h)v}^2+\frac{3}{2}\mynormb{v}_{-1}^2$.
\end{lemma}
\begin{proof}
The generalized H\"{o}lder inequality \eqref{ieq:generalized Holder} and the Young's inequality lead to
\[
\mynormb{v}^2\le \mynormb{\nabla_hv}\mynormb{v}_{-1}
\le \frac{\varepsilon_1}{2}\mynormb{\nabla_hv}^2
+ \frac{1}{2\varepsilon_1}\mynormb{v}_{-1}^2\quad\text{for $\varepsilon_1>0$}.
\]
Also, by using the discrete Green's formula and Cauchy-Schwarz inequality, one has
\begin{align*}
\mynormb{\nabla_hv}^2
=\mynormb{v}^2 - \myinnerb{(1+\Delta_h)v,v}\le \bra{1+\frac{\varepsilon_2}{2}}\mynormb{v}^2
+\frac{1}{2\varepsilon_2}\mynormb{(1+\Delta_h)v}^2
\quad\text{for $\varepsilon_2>0$}.
\end{align*}
The above two inequalities with
$\varepsilon_1=\frac{2}{3}$ and $\varepsilon_2=1$ yields the claimed result.
\end{proof}

Note that, the solution $\phi^n$ of BDF2 scheme \eqref{scheme:PFC BDF2}
preserves the volume, $\myinnerb{\phi^n,1}=\myinnerb{\phi^0,1}$,
for $n\ge1$. Actually, taking the inner product of \eqref{scheme:PFC BDF2} by 1 and applying the summation by parts, one has
 $\myinnerb{D_2\phi^j,1}=\myinnerb{\Delta_h\mu^j,1}=0$ for $j\ge1$.
 Multiplying both sides of this equality by the DOC kernels
 $\theta_{n-j}^{(n)}$ and summing the index $j$ from $j=1$ to $n$, we get
$$\sum_{j=1}^n\theta_{n-j}^{(n)}\myinnerb{D_2\phi^j,1}=0\quad\text{for $n\ge1$}.$$
It leads to $\myinnerb{\diff \phi^n,1}=0$ directly by taking $v^j=\phi^j$ in the equality \eqref{orthogonal equality for BDF2}.
Simple induction yields the volume conversation law, $\myinnerb{\phi^n,1}=\myinnerb{\phi^{n-1},1}=\cdots=\myinnerb{\phi^{0},1}$
for $n\ge1$.

\begin{theorem}\label{thm: convexity solvability}
If the step size $\tau_n\le\frac{2+4r_n}{3\epsilon(1+r_n)}$,
the BDF2 scheme \eqref{scheme:PFC BDF2} is uniquely solvable.
\end{theorem}

\begin{proof}
For any fixed time-level indexes $n\ge1$,
we consider the following energy functional $G$ on the space
$\mathbb{V}_{h}^{*}:=\big\{z\in\mathbb{V}_{h}\,|\, \myinnerb{z,1}=\myinnerb{\phi^{n-1},1}\big\},$
\begin{align*}
G[z]:=\frac{1}{2}b_0^{(n)}\mynormb{z-\phi^{n-1}}_{-1}^2+b_1^{(n)}\myinnerb{\diff \phi^{n-1},z}_{-1}
+\frac12\mynormb{(1+\Delta_h)z}^2
+\frac14\mynormb{z}_{l^4}^4-\frac{\epsilon}2\mynormb{z}^2.
\end{align*}
Under the time-step size condition $\tau_n\le\frac{2+4r_n}{3\epsilon(1+r_n)}$
or $b_0^{(n)}\ge3\epsilon/2$,
the functional $G$ is strictly convex since, for any $\lambda\in \mathbb{R}$ and any $\psi\in \mathbb{\mathring V}_{h}$,
\begin{align*}
\frac{\zd^2G}{\zd\lambda^2}[z+\lambda\psi]\Big|_{\lambda=0}=&\,b_0^{(n)}\mynormb{\psi}_{-1}^2
+\mynormb{(1+\Delta_h)\psi}^2
+3\mynormb{z\psi}^2-\epsilon\mynormb{\psi}^2\\
\ge&\,\brab{b_0^{(n)}-\frac{3\epsilon}2}\mynormb{\psi}_{-1}^2
+\frac{2}{3}\mynormb{(1+\Delta_h)\psi}^2
+3\mynormb{z\psi}^2>0,
\end{align*}
where Lemma \ref{lem: L2-Embedding-Inequality} has been applied with the setting $0<\epsilon<1$.
Thus the functional $G$ has a unique minimizer, denoted by $\phi^n$, if and only if it solves the equation
\begin{align*}
0=\frac{\zd G}{\zd\lambda}[z+\lambda\psi]\Big|_{\lambda=0}=&\,
\myinnerb{b_0^{(n)}\brat{z-\phi^{n-1}}+b_1^{(n)}\diff \phi^{n-1},\psi}_{-1}
+\myinnerb{(1+\Delta_h)^2z+z^3-\epsilon z,\psi}\\
=&\,
\myinnerB{b_0^{(n)}\brat{z-\phi^{n-1}}+b_1^{(n)}\diff \phi^{n-1}
-\Delta_h\bra{(1+\Delta_h)^2z+z^3-\epsilon z},\psi}_{-1}.
\end{align*}
This equation holds for any $\psi\in \mathbb{\mathring V}_{h}$ if and only
if the unique minimizer $\phi^n\in\mathbb{V}_{h}^{*}$ solves
\begin{align*}
b_0^{(n)}\brat{\phi^n-\phi^{n-1}}+b_1^{(n)}\diff \phi^{n-1}
-\Delta_h\bra{(1+\Delta_h)^2\phi^n+(\phi^n)^3-\epsilon \phi^n}=0,
\end{align*}
which is just the BDF2 scheme \eqref{scheme:PFC BDF2}.
It verifies the claimed result  and completes the proof.
\end{proof}

The proof of Theorem \ref{thm: convexity solvability} also says that
the BDF2 scheme \eqref{scheme:PFC BDF2} is equivalent to the minimization
of a convex functional $G[z]$ under the condition $\tau_n\le \frac{2+4r_n}{3\epsilon(1+r_n)}$.
We see that the BDF2 implicit time-stepping scheme is also convex according to
Xu et al. \cite{xu2019on}.

\subsection{Energy dissipation law}
The following result, cf. \cite[Lemma 2.1]{Liao2019Adaptive},
shows that the BDF2 convolution kernels $b_{n-k}^{(n)}$ are positive definite
provided the adjacent time-step ratios $r_k$ satisfy \textbf{S1},
or $0<r_k< r_{\mathrm{sup}}$, where $r_{\mathrm{sup}}=\frac{3+\sqrt{17}}{2}$
is the positive root of the equation $2+3r_{\mathrm{sup}}-r_{\mathrm{sup}}^2=0$.
Consider the function
\begin{align}\label{def: step-ratios function}
R(z,s):=\frac{2+4z-z^2}{1+z}-\frac{s}{1+s}\quad
\text{for $0\le z,s< r_{\mathrm{sup}}.$}
\end{align}
It is easy to check that $R(z,s)$ is increasing in $(0,\sqrt{3}-1)$
and decreasing in $(\sqrt{3}-1, r_{\mathrm{sup}})$ with respect to $z$,
and decreasing with respect to $s$.
So the condition \textbf{S1} ensures
\begin{align*}
\frac{2+4r_k-r_k^2}{1+r_k} - \frac{r_{k+1}}{1+r_{k+1}}
=R(r_k,r_{k+1})>R(r_k,r_{\mathrm{sup}})>0\quad \text{for $k\ge1$.}
\end{align*}


\begin{lemma}\label{lem:conv kernels positive}
Let \textbf{S1} holds. For any real sequence $\{w_k\}_{k=1}^n$ with n entries,
it holds that
\begin{align*}
2w_k\sum_{j=1}^kb_{k-j}^{(k)}w_j
&\ge\frac{r_{k+1}}{1+r_{k+1}}\frac{w_k^2}{\tau_k}
-\frac{r_k}{1+r_k}\frac{w_{k-1}^2}{\tau_{k-1}}
+R(r_k,r_{k+1})
\frac{w_k^2}{\tau_k}
\end{align*}
for $k\ge2$.
So the discrete convolution kernels $b_{n-k}^{(n)}$ are positive definite,
\[
\sum_{k=1}^n w_k \sum_{j=1}^k b_{k-j}^{(k)}w_j\ge
\frac12\sum_{k=1}^n R(r_k,r_{k+1})
\frac{w_k^2}{\tau_k}\quad\text{for $n\ge 1$}.
\]
\end{lemma}

Now we prove the energy stability of BDF2 scheme \eqref{scheme:PFC BDF2}.
Let $E[\phi^k]$ be the discrete version of free energy functional \eqref{cont:free energy}, given by
\begin{align}\label{def: discrete free energy}
E[\phi^k]
:=\frac{1}{2}\mynormb{(1+\Delta_h)\phi^k}^2 +
\frac{1}{4}\mynormb{(\phi^k)^2-\epsilon}^2-\frac{1}{4}\mynormb{\epsilon}^2 \quad\text{for $k\ge 0$.}
\end{align}
Since the BDF2 formula \eqref{def: BDF2-Formula} is naturally self-dissipative, we define a modified discrete energy,
\begin{align*}
\mathcal{E}[\phi^k]
:=E[\phi^k] + \frac{r_{k+1}}{2(1+r_{k+1})\tau_k}\mynormb{\diff \phi^k}_{-1}^2
 \quad\text{for $k\ge 0$}
\end{align*}
where $\mathcal{E}[\phi^0]=E[\phi^0]$ due to the setting $r_1\equiv0$.


\begin{theorem}\label{thm:energy decay law}
Assume that \textbf{S1} holds and the time-step sizes are properly small such that
\begin{align}\label{Restriction-Time-Step}
\tau_n \le\frac{2}{3\epsilon}\min\Big\{\frac{1+2r_n}{1+r_n},R(r_n,r_{n+1})
\Big\}
\quad\text{for $n\ge 1$,}
\end{align}
the variable-step BDF2 scheme \eqref{scheme:PFC BDF2} preserves the following energy dissipation law
\begin{align*}
\mathcal{E}[\phi^n] \le \mathcal{E}[\phi^{n-1}]\le \mathcal{E}[\phi^{0}]=E[\phi^{0}]
\quad\text{for $n\ge 1$.}
\end{align*}
\end{theorem}
\begin{proof}
The first condition of \eqref{Restriction-Time-Step} ensures the unique solvability in Theorem \ref{thm: convexity solvability}.
We will establish the energy dissipation law under the second condition of \eqref{Restriction-Time-Step}.
The volume conversation law implies $\diff \phi^n\in\mathbb{\mathring V}_{h}$ for $n\ge1$.
Then we make the inner product of \eqref{scheme:PFC BDF2} by $(-\Delta_h)^{-1}\diff \phi^{n}$ and obtain
\begin{align}\label{Energy-Law-Inner}
\myinnerb{D_2\phi^n,(-\Delta_h)^{-1}\diff \phi^{n}}
+\myinnerb{(1+\Delta_h)^2\phi^n,\diff \phi^{n}}
+\myinnerb{\brat{\phi^n}^3-\epsilon\phi^n,\diff \phi^{n}}=0.
\end{align}
With the help of the summation by parts and $2a(a-b)=a^2-b^2+(a-b)^2$,
the second term at the left hand side of \eqref{Energy-Law-Inner} gives
\begin{align*}
\myinnerb{(1+\Delta_h)^2\phi^n,\diff \phi^n}
=\frac{1}{2}\mynormb{(1+\Delta_h)\phi^n}^2
-\frac{1}{2}\mynormb{(1+\Delta_h)\phi^{n-1}}^2
+\frac{1}{2}\mynormb{(1+\Delta_h)\diff \phi^{n}}^2.
\end{align*}
It is easy to check the following identity
\begin{align*}
4\brab{a^3-\epsilon a}\bra{a-b}= \bra{a^2-\epsilon}^2 - \bra{b^2-\epsilon}^2-2\bra{\epsilon-a^2}\bra{a-b}^2 + \bra{a^2-b^2}^2.
\end{align*}
Then the third term in \eqref{Energy-Law-Inner} can be bounded by
\begin{align*}
\myinnerb{(\phi^n)^3-\epsilon\phi^n,\diff \phi^n}
&\ge \frac{1}{4}\mynormb{(\phi^n)^2-\epsilon}^2 - \frac{1}{4}\mynormb{(\phi^{n-1})^2-\epsilon}^2
-\frac{1}{2}\myinnerb{(\epsilon-(\phi^n)^2)\bra{\diff \phi^n}^2,1}\\
&\ge \frac{1}{4}\mynormb{(\phi^n)^2-\epsilon}^2 - \frac{1}{4}\mynormb{(\phi^{n-1})^2-\epsilon}^2
-\frac{\epsilon}{2}\mynormb{\diff \phi^n}^2.
\end{align*}
Thus it follows from \eqref{Energy-Law-Inner} that
\begin{align*}
\myinnerb{D_2\phi^n,(-\Delta_h)^{-1}\diff \phi^{n}}
+\frac{1}{2}\mynormb{(1+\Delta_h)\diff \phi^{n}}^2
-\frac{\epsilon}{2}\mynormb{\diff \phi^n}^2+E[\phi^n]\le E[\phi^{n-1}].
\end{align*}
Applying Lemma \ref{lem: L2-Embedding-Inequality}, one has
$$\frac{\epsilon}{2}\mynormb{\diff \phi^n}^2 \le \frac{1}{6}\mynormb{(1+\Delta_h)\diff \phi^n}^2
+\frac{3\epsilon}{4}\mynormb{\diff \phi^n}_{-1}^2,
$$
where $0<\epsilon<1$ has been used. Thus we can obtain that
\begin{align}\label{Energy-Inequality}
\myinnerb{D_2\phi^n,(-\Delta_h)^{-1}\diff \phi^{n}}
-\frac{3\epsilon}{4}\mynormb{\diff \phi^n}_{-1}^2
+E[\phi^n]\le E[\phi^{n-1}]\quad\text{for $n\ge 1$.}
\end{align}

For the general cases $n\ge 2$, we take $w_j=\diff \phi^{j}$
in the first inequality of Lemma \ref{lem:conv kernels positive}
and apply the condition $\frac{1}{2\tau_n}R(r_n,r_{n+1})\ge \frac{3}{4}\epsilon$
to obtain
\begin{align*}
\myinnerb{D_2\phi^n,(-\Delta_h)^{-1}\diff \phi^{n}}
\ge&\, \frac{r_{n+1}\tau_n}{2(1+r_{n+1})}\mynormb{\partial_\tau\phi^n}_{-1}^2
-\frac{r_n\tau_{n-1}}{2(1+r_n)}\mynormb{\partial_\tau\phi^{n-1}}_{-1}^2
+\frac{3\epsilon}{4}\mynormb{\diff \phi^n}_{-1}^2.
\end{align*}
Combining the above inequality with \eqref{Energy-Inequality} yields
$\mathcal{E}[\phi^n] \le \mathcal{E}[\phi^{n-1}]$ for $n\ge 2$.
It remains to consider $n=1$.
Recalling $b_0^{(1)}=1/\tau_1$, we use
$\frac{1}{2\tau_1}R(r_1,r_2)=\frac{2+r_{2}}{2\tau_1(1+r_{2})}\ge \frac{3}{4}\epsilon$
to derive that
\begin{align*}
\myinnerb{D_2\phi^1,(-\Delta_h)^{-1}\diff \phi^{1}}
=&\,\myinnerb{D_1\phi^1,(-\Delta_h)^{-1}\diff \phi^{1}}
=\braB{\frac{r_{2}\tau_1}{2(1+r_2)}
+\frac{(2+r_{2})\tau_1}{2(1+r_2)}}\mynormb{\partial_\tau\phi^1}_{-1}^2\\
\ge&\,\frac{r_{2}\tau_1}{2(1+r_2)}\mynormb{\partial_\tau\phi^1}_{-1}^2
+\frac{3\epsilon}{4}\mynormb{\diff \phi^1}_{-1}^2.
\end{align*}
We combine this inequality with \eqref{Energy-Inequality} to find
$\mathcal{E}[\phi^1]\le\mathcal{E}[\phi^0]=E[\phi^0]$, and complete the proof.
\end{proof}

\begin{remark}\label{remark:comments on time-step condition}
Some remarks on the time-step size constraint
\eqref{Restriction-Time-Step} are listed here under the step-ratio condition \textbf{S1},
that is, $0<r_k< r_{\mathrm{sup}}$ for $k\ge2$.
For $n=1$, it gives $\epsilon\tau_1\le \frac2{3}\min\{1,\frac{2+r_2}{1+r_2}\}=\frac23$
and one can choose $\tau_n$ such that $\epsilon\tau_1\le\frac{2}3$.
Recalling the monotonicity of $R(z,s)$, we consider the following three cases for $n\ge2$:
\begin{itemize}
  \item[(i)] If $0<r_n,r_{n+1}\le \sqrt{3}-1$, $R(r_n,r_{n+1})\ge R(0,\sqrt{3}-1)$
  and then $\epsilon\tau_n\le\frac{2}{3}
  \min\big\{1,\frac{3+\sqrt{3}}{3}\big\}=\frac{2}{3}$.
  One can choose the step size $\tau_n\le\frac{2}{3\epsilon}$;
  \item[(ii)] If $\sqrt{3}-1<r_n,r_{n+1}\le 2$, $R(r_n,r_{n+1})\ge R(2,2)$ and then $\epsilon\tau_n\le\frac{2}{3}\min\big\{\frac{6-\sqrt{3}}{3},\frac{4}{3}\big\}=\frac{8}{9}$.
  One can choose the step size  $\tau_n\le\frac{8}{9\epsilon}$;
  \item[(iii)] If $2<r_n< r_{\mathrm{sup}}$, one can choose a small step size $\tau_{n+1}$ or step ratio $r_{n+1}$ to ensure the step size restriction \eqref{Restriction-Time-Step}
  in adaptive computations, especially when the current step-ratio $r_n\rightarrow r_{\mathrm{sup}}$.
   For an example, the choice $\tau_n\le\frac{1}{4\epsilon}$ is sufficient
  if one choose $R(r_{\mathrm{sup}},r_{n+1})\ge\frac3{8}$, i.e., the next time-step ratio  $r_{n+1}\le(3+16\sqrt{17})/101\approx0.68$.
\end{itemize}
In summary, the time-step size constraint
\eqref{Restriction-Time-Step}
is always mild in practical computations.
\end{remark}

\begin{lemma}\label{lem:Bound-Solution}
Assume that \textbf{S1} holds and the time--step sizes fulfill
\eqref{Restriction-Time-Step}.
The solution of BDF2 scheme \eqref{scheme:PFC BDF2} is stable in the $L^{\infty}$ norm,
$$\mynormb{\phi^n}_{\infty}\le c_0:=c_{\Omega}\sqrt{8E[\phi^{0}]+2\brab{2+\epsilon}^2\abs{\Omega_h}}\quad\text{for $n\ge 1$},$$
where $c_0$ is dependent on the domain $\Omega$ and the initial value $\phi^0$,
but independent of the time $t_n$, step sizes $\tau_n$ and step ratios $r_n$.
\end{lemma}

\begin{proof}
Since $(a^2-2-\epsilon)^2\ge0$,  one has
$\mynorm{\phi^n}_{l^4}^4\ge\brab{4+2\epsilon}\mynormb{\phi^n}^2
-\brab{2+\epsilon}^2\absb{\Omega_h}$.
The energy dissipation law in Theorem \ref{thm:energy decay law} shows that
$E[\phi^{0}]\ge \mathcal{E}[\phi^{n}]\ge E[\phi^{n}]$. Then we have
\begin{align*}
8E[\phi^{0}]\ge 8E[\phi^{n}]
=&\, 4\mynormb{(1+\Delta_h)\phi^n}^2
+2\mynormb{\phi^n}_{l^4}^4-4\epsilon\mynormb{\phi^n}^2\\
\ge&\, 4\mynormb{(1+\Delta_h)\phi^n}^2+8\mynormb{\phi^n}^2
-2\brab{2+\epsilon}^2\absb{\Omega_h}\\
\ge&\, 2\mynormb{\Delta_h\phi^n}^2+4\mynormb{\phi^n}^2
-2\brab{2+\epsilon}^2\absb{\Omega_h}\\
\ge&\,\bra{\mynormb{\Delta_h\phi^n}+\mynormb{\phi^n}}^2
-2\brab{2+\epsilon}^2\absb{\Omega_h}
\quad\text{for $n\ge1$,}
\end{align*}
where the inequality, $\mynormb{\Delta_h\phi^n}^2
=\mynormb{(1+\Delta_h)\phi^n-\phi^n}^2
\le 2\mynormb{(1+\Delta_h)\phi^n}^2+2\mynormb{\phi^n}^2$,
was used in the second inequality.
Then one can apply the Sobolev embedding inequality \eqref{ieq: Linfty H2 embedding}
to obtain
\begin{align*}
\mynormb{\phi^n}_{\infty}^2
\le c_{\Omega}^2\brab{\mynormb{\phi^n}+\mynormb{\Delta_h\phi^n}}^2
\le c_{\Omega}^2\brab{8E[\phi^{0}]+2\brab{2+\epsilon}^2\abs{\Omega_h}}=c_0^2,
\end{align*}
which leads to the claimed bound and completes the proof.
\end{proof}

\section{$L^2$ norm error estimate}
\setcounter{equation}{0}

\subsection{Some properties of DOC kernels}

The following lemma gathers the results of Lemma 2.2, Corollary 2.1 and Lemma 2.3 in \cite{Liao2019Adaptive}.
\begin{lemma}\label{lem:DOC property}
	If the discrete convolution kernels  $b^{(n)}_{n-k}$ defined in
	\eqref{def:BDF2-kernels} are positive semi-definite (the restriction \textbf{S1} is sufficient),
then the DOC kernels $\theta_{n-j}^{(n)}$
	defined in \eqref{DOC-Kernels} satisfy:
\begin{itemize}
  \item[(I)] The discrete kernels $\theta_{n-j}^{(n)}$ are positive definite;
  \item[(II)] The discrete kernels $\theta_{n-j}^{(n)}$ are positive and
  \begin{align*}
	\theta_{n-j}^{(n)}=\frac{1}{b^{(j)}_{0}}\prod_{i=j+1}^n\frac{r_i^2}{1+2r_i}
	\quad\text{ for $1\le j\le n$;}
	\end{align*}
  \item[(III)] $\displaystyle \sum_{j=1}^{n}\theta_{n-j}^{(n)}=\tau_n$ such that
  $\displaystyle \sum_{k=1}^{n}\sum_{j=1}^{k}\theta_{k-j}^{(k)}=t_n$ for $n\ge1$.
\end{itemize}
\end{lemma}

To facilitate the convergence analysis, we present a discrete convolution
inequality with respect to the DOC kernels $\theta_{n-j}^{(n)}$,
but leave the proof to Appendix A.

\begin{lemma}\label{lem:quadr form inequ}
If \textbf{S1} holds, then for any real sequences $\{w_k\}_{k=1}^n$
and $\{v_k\}_{k=1}^n$,
\begin{align*}
\sum_{k=1}^n  \sum_{j=1}^k\theta_{k-j}^{(k)} w_k v_j
\le\varepsilon\sum_{k=1}^n  \sum_{j=1}^k\theta_{k-j}^{(k)} v_k v_j
+\frac{\mathcal{M}_r}{\varepsilon}\sum_{k=1}^n  \sum_{j=1}^k\theta_{k-j}^{(k)} w_k w_j
\quad \text{$\forall\;\varepsilon > 0$},
\end{align*}
where  $\mathcal{M}_r>0$ is a constant independent of the time $t_n$,
time-step sizes $\tau_n$ and step ratios $r_n$.
\end{lemma}


Lemma \ref{lem:quadr form inequ} yields
the following discrete embedding-like inequality in the quadratic form. Here and hereafter,
we use the notation
$\sum_{k,j}^{n,k}\triangleq\sum_{k=1}^n\sum_{j=1}^k$ for the sake of brevity.

\begin{lemma}\label{lem:inner quad ineq}
If \textbf{S1} holds, then for any grid function $v^n\in\mathbb{V}_{h}$
and any constant $\varepsilon>0$,
\begin{align*}
\sum_{k,j}^{n,k}\theta_{k-j}^{(k)}\myinnerb{\Delta_hv^j,\Delta_h v^k}
\le&\,\frac{16\mathcal{M}_r^3}{\varepsilon^2}
\sum_{k,j}^{n,k}\theta_{k-j}^{(k)}\myinnerb{v^j,v^k}
+\varepsilon\sum_{k,j}^{n,k}\theta_{k-j}^{(k)}
\myinnerb{\nabla_h\Delta_hv^j,\nabla_h\Delta_h v^k}.
\end{align*}
\end{lemma}

\begin{proof}For any constant $\varepsilon_3>0$, we can take
$w_k:=-\nabla_h\Delta_h v^k$, $v_j:=-\nabla_hv^j$ and $\varepsilon:=\mathcal{M}_r/\varepsilon_3$ in
Lemma \ref{lem:quadr form inequ} and derive that
\begin{align*}
2\sum_{k,j}^{n,k}\theta_{k-j}^{(k)}\myinnerb{\Delta_hv^j,\Delta_h v^k}
\le&\,\frac{2\mathcal{M}_r}{\varepsilon_3}\sum_{k,j}^{n,k}\theta_{k-j}^{(k)}\myinnerb{\nabla_hv^j,\nabla_hv^k}
+2\varepsilon_3\sum_{k,j}^{n,k}\theta_{k-j}^{(k)}
\myinnerb{\nabla_h\Delta_hv^j,\nabla_h\Delta_h v^k}
\end{align*}
Similarly, Lemma \ref{lem:quadr form inequ}
with  $v_j:=-\Delta_h v^j$, $w_k:=v^k$  and $\varepsilon:=\varepsilon_3/(2\mathcal{M}_r)$ yields
\begin{align*}
\frac{2\mathcal{M}_r}{\varepsilon_3}\sum_{k,j}^{n,k}\theta_{k-j}^{(k)}\myinnerb{\nabla_hv^j,\nabla_hv^k}
\le&\,
\sum_{k,j}^{n,k}\theta_{k-j}^{(k)}
\myinnerb{\Delta_hv^j,\Delta_h v^k}
+\frac{4\mathcal{M}_r^3}{\varepsilon_3^2}\sum_{k,j}^{n,k}\theta_{k-j}^{(k)}\myinnerb{v^j,v^k}.
\end{align*}
We complete the proof by summing up the above two inequalities
and taking $\varepsilon_3=\varepsilon/2$.
\end{proof}

\subsection{Convolutional consistency in time}

Now consider the error behavior of BDF2 time-stepping
with respect to the variation of time-step sizes.
Let $\xi_{\Phi}^j:=D_2\Phi(t_j)-\partial_t\Phi(t_j)$ be the local consistency error of BDF2 formula at the time $t=t_j$.
We will consider a convolutional consistency error $\Xi_{\Phi}^k$ defined by
\begin{align}\label{BDF2-global consistency}
\Xi_{\Phi}^k:=\sum_{j=1}^k\theta_{k-j}^{(k)}\xi_{\Phi}^j=
\sum_{j=1}^k\theta_{k-j}^{(k)}\bra{D_2\Phi(t_j)-\partial_t\Phi(t_j)}\quad\text{for $k\ge1$.}
\end{align}
\begin{lemma}\label{lem:BDF2-Consistency-Error}
If \textbf{S1} holds, the convolutional consistency error $\Xi_{\Phi}^k$ in \eqref{BDF2-global consistency} satisfies
\begin{align*}
\absb{\Xi_{\Phi}^k}
\le&\, \theta_{k-1}^{(k)}\int_{0}^{t_1} \absb{\Phi''(s)} \zd{s}
+3\sum_{j=1}^{k}\theta_{k-j}^{(k)}\tau_{j}\int_{t_{j-1}}^{t_j} \absb{\Phi'''(s)} \zd{s}\quad\text{for $k\ge1$}
\end{align*}
such that
\begin{align*}
\sum_{k=1}^n\absb{\Xi_{\Phi}^k}
\le&\, \tau_1\int_{0}^{t_1} \absb{\Phi''(s)} \zd{s}\,\sum_{k=1}^n\prod_{i=2}^k\frac{r_i^2}{1+2r_i}
+3t_n\max_{1\le j\le n}\braB{\tau_{j}\int_{t_{j-1}}^{t_j} \absb{\Phi'''(s)} \zd{s}}\quad\text{for $n\ge1$.}
\end{align*}
\end{lemma}
\begin{proof}We use the notations $G_{t2}^{j}=\int_{t_{j-1}}^{t_j}\absb{\Phi''(s)}\zd{s}$ and
$G_{t3}^{j}=\int_{t_{j-1}}^{t_j}\absb{\Phi'''(s)}\zd{s}$ for $j\ge1$.
The proof of \cite[Lemma 3.2]{Liao2019Adaptive} gives
\begin{align*}
\absb{\xi_{\Phi}^1}\le b_0^{(1)}\tau_1G_{t2}^{1}\quad\text{and}\quad
\absb{\xi_{\Phi}^j}\le b_0^{(j)}\tau_j^2G_{t3}^{j}+\frac{r_j^2}{2(1+2r_j)}b_0^{(j)}\tau_{j-1}^2G_{t3}^{j-1}\quad\text{for $j\ge2$}.
\end{align*}
Recalling the definitions of BDF2 kernels \eqref{def:BDF2-kernels} and DOC kernels \eqref{DOC-Kernels}, we find that
\begin{align*}
\theta_{k-j}^{(k)}b_0^{(j)}=-\theta_{k-j-1}^{(k)}b_1^{(j+1)}=
 \frac{r_{j+1}^2}{1+2r_{j+1}}\theta_{k-j-1}^{(k)}b_0^{(j+1)}\quad\text{for $1\le j\le k-1$}.
\end{align*}
Lemma \ref{lem:DOC property}(II) shows $\theta_{k-j}^{(k)}>0$. Thus we apply the triangle inequality to derive that
\begin{align*}
\absb{\Xi_{\Phi}^k}\le&\,\sum_{j=1}^k\theta_{k-j}^{(k)}\absb{\xi_{\Phi}^k}
=\theta_{k-1}^{(k)}\absb{\xi_{\Phi}^1}+\sum_{j=2}^k\theta_{k-j}^{(k)}\absb{\xi_{\Phi}^k}\\
\le&\,\theta_{k-1}^{(k)}b_0^{(1)}\tau_1G_{t2}^{1}+
\sum_{j=2}^k\theta_{k-j}^{(k)}b_0^{(j)}\tau_j^2G_{t3}^{j}
+\frac12\sum_{j=1}^{k-1}\theta_{k-j-1}^{(k)}b_0^{(j+1)}\frac{r_{j+1}^2}{1+2r_{j+1}}\tau_{j}^2G_{t3}^{j}\\
=&\,\theta_{k-1}^{(k)}b_0^{(1)}\tau_1G_{t2}^{1}+
\sum_{j=2}^k\theta_{k-j}^{(k)}b_0^{(j)}\tau_j^2G_{t3}^{j}
+\frac12\sum_{j=1}^{k-1}\theta_{k-j}^{(k)}b_0^{(j)}\tau_{j}^2G_{t3}^{j}\\
\le&\,\theta_{k-1}^{(k)}b_0^{(1)}\tau_1G_{t2}^{1}
+\frac32\sum_{j=1}^{k}\theta_{k-j}^{(k)}b_0^{(j)}\tau_{j}^2G_{t3}^{j}\quad\text{for $k\ge1$}.
\end{align*}
Moreover, the definition \eqref{def:BDF2-kernels} of BDF2 kernels yields that $b_0^{(1)}\tau_1=1$ and $b_0^{(j)}\tau_{j}=\frac{1+2r_j}{1+r_j}\le2$,
and the claimed first inequality follows immediately.
Then the second estimation can be derived by using Lemma \ref{lem:DOC property} (II)-(III). It completes the proof.
\end{proof}

\begin{remark}\label{remark: S1-S2 accuracy order}
Under the mild condition \textbf{S1}, Lemma \ref{lem:DOC property} (III) shows $\theta_{k-1}^{(k)}\le \tau_k$ such that
\begin{align*}
\sum_{k=1}^n\absb{\Xi_{\Phi}^k}\le&\, t_n\int_{0}^{t_1} \absb{\Phi''(s)}\zd{s}
+3t_n\max_{1\le j\le n}\braB{\tau_{j}\int_{t_{j-1}}^{t_j} \absb{\Phi'''(s)} \zd{s}}\quad\text{for $n\ge1$.}
\end{align*}
It will arrive at first-order convergence. The degraded accuracy is mainly attributed to the convolutional accumulation,
from the irregular variation of time-step sizes, onto the first-level solution.
If \textbf{S2} holds, there exists a bounded quantity $c_r=c_r(N_0,r_{c},\hat{r}_{c})$ such that
\begin{align*}
\sum_{k=1}^n\prod_{i=2}^k\frac{r_i^2}{1+2r_i}
\le c_r(N_0,r_{c},\hat{r}_{c}):=\braB{\frac{\hat{r}_{c}^2}{1+2\hat{r}_{c}}}^{N_0}\frac{1+2r_{c}}{1+2r_{c}-r_{c}^2}\,,
\end{align*}
where $r_c$ takes the maximum value of all step ratios $r_k\in (0,1+\sqrt{2})$
and $\hat{r}_c$ takes the maximum value of those step ratios $r_k\in\big[1+\sqrt{2},\frac{3+\sqrt{17}}{2}\big)$ for $2\le k\le N$. Then Lemma \ref{lem:BDF2-Consistency-Error} shows
\begin{align*}
\sum_{k=1}^n\absb{\Xi_{\Phi}^k}\le&\, c_r(N_0,r_{c},\hat{r}_{c})\tau_1\int_{0}^{t_1} \absb{\Phi''(s)}\zd{s}
+3t_n\max_{1\le j\le n}\braB{\tau_{j}\int_{t_{j-1}}^{t_j} \absb{\Phi'''(s)} \zd{s}}\quad\text{for $n\ge1$,}
\end{align*}
which will yield the desired second-order convergence.
\end{remark}

\subsection{Convergence analysis}
We use  the standard semi-norms and norms in the Sobolev space $H^{m}(\Omega)$ for $m\ge0$.
Let $C_{per}^{\infty}(\Omega)$ be a set of infinitely differentiable $L$-periodic functions defined on $\Omega$,
and $H_{per}^{m}(\Omega)$ be the closure of $C_{per}^{\infty}(\Omega)$ in $H^{m}(\Omega)$,
endowed with the semi-norm $|\cdot|_{H_{per}^m}$ and the norm $\mynorm{\cdot}_{H_{per}^{m}}$.

For simplicity, denote $|\cdot|_{H^m}:=|\cdot|_{H_{per}^m}$, $\mynorm{\cdot}_{H^{m}}:=\mynorm{\cdot}_{H_{per}^{m}}$, and $\mynorm{\cdot}_{L^{2}}:=\mynorm{\cdot}_{H^{0}}$.
We denote the maximum norm by $\mynorm{\cdot}_{L^{\infty}}$
and have the Sobolev embedding inequality $\mynorm{u}_{L^{\infty}}\le C_{\Omega}\mynorm{u}_{H^{2}}$
for $u\in C_{per}^{\infty}(\Omega)\cap H_{per}^{m}(\Omega)$.
Next lemma lists some approximations, see \cite{Shen2011Spectral}, of the $L^2$-projection operator $P_{M}$ and
trigonometric interpolation operator $I_{M}$ defined in subsection 2.1.
\begin{lemma}\label{lem:Projection-Estimate}
For any $u\in{H_{per}^{q}}(\Omega)$ and $0\le{s}\le{q}$, it holds that
\begin{align}
\mynorm{P_{M}u-u}_{H^{s}}
\le C_uh^{q-s}|u|_{H^{q}},
\quad \mynorm{P_{M}u}_{H^{s}}\le C_u\mynorm{u}_{H^{s}};\label{Projection-Estimate}
\end{align}
and, in addition if $q>3/2$,
\begin{align}
\mynorm{I_{M}u-u}_{H^{s}}
\le C_uh^{q-s}|u|_{H^{q}},
\quad \mynorm{I_{M}u}_{H^{s}}\le C_u\mynorm{u}_{H^{s}}.\label{Interpolation-Estimate}
\end{align}
\end{lemma}

Note that,  the energy dissipation law \eqref{cont:energy dissipation}
of PFC model \eqref{cont: Problem-PFC} shows that $E[\Phi^n]\le E[\Phi(t_0)]$.
From the formulation \eqref{cont:free energy}, it is not difficult to see that
$\mynormb{\Phi^n}_{H^2}$ can be bounded by a time-independent constant.
By the projection estimate \eqref{Projection-Estimate} in Lemma \ref{lem:Projection-Estimate}
and the Sobolev inequality, one has
$\mynormb{P_M\Phi^n}_{L^{\infty}}\le c_1$ and then
\begin{align}\label{maximum bound projection}
\mynormb{P_M\Phi^n}_{\infty}\le\mynormb{P_M\Phi^n}_{L^{\infty}}\le c_1\quad\text{for $1\le n\le N$,}
\end{align}
where $c_1$ is dependent on the domain $\Omega$
and initial data $\Phi(t_0)$, but independent of the time $t_n$.
We are in the position to prove the $L^2$ norm convergence of the adaptive BDF2 scheme \eqref{scheme:PFC BDF2}
by choosing an initial value $\phi^0=I_M\Phi(t_0)$.

\begin{theorem}\label{thm:Convergence-Results}
Assume that the PFC problem \eqref{cont: Problem-PFC} has a solution $\Phi\in C^3\brab{[0,T];{H}_{per}^{m+6}}$
for some integer $m\ge0$. Suppose further that the step-ratios condition \textbf{S1} and
the time-step size restriction \eqref{Restriction-Time-Step} hold such that
the adaptive BDF2 implicit scheme \eqref{scheme:PFC BDF2} is unique solvable and energy stable.
If the maximum time-step size $\tau$ is sufficiently small such that $\tau\le 1/\bra{2c_2}$,
then the solution $\phi^n$ is (at least, first-order) convergent in the $L^2$ norm,
\begin{align*}
\mynorm{\Phi^n-\phi^n}
\le C_\phi\exp\bra{2c_2t_{n-1}}&\,\bigg[(1+t_n)h^m
+\tau_1\int_{0}^{t_1} \mynormb{\Phi''(s)}_{L^2}\zd{s}\,\sum_{k=1}^n\prod_{i=2}^k\frac{r_i^2}{1+2r_i}\\
&\,+3t_n\max_{1\le j\le n}\braB{\tau_{j}\int_{t_{j-1}}^{t_j} \mynormb{\Phi'''(s)}_{L^2} \zd{s}}\bigg]\quad\text{for $1\le n\le N$,}
\end{align*}
where the positive constant $c_2:=250\mathcal{M}_r^3+4\mathcal{M}_r(c_1^2+c_0c_1+c_0^2+\epsilon)^2$
is always dependent on the domain $\Omega$ and
the initial values $\Phi^0$ and $\phi^0$,
but independent of the time $t_n$, step sizes $\tau_n$ and step ratios $r_n$.
Remark \ref{remark: S1-S2 accuracy order} shows that the second-order time accuracy will be recovered
by replacing the weak step-ratio condition \textbf{S1} with a mild one \textbf{S2}.
\end{theorem}

\begin{proof}
We evaluate the $L^2$ norm error $\mynorm{\Phi^n-\phi^n}$ by a usual splitting,
$$\Phi^n-\phi^n=\Phi^n-\Phi_M^n+e^n,$$
where $\Phi_M^n:=P_M\Phi^n$ is the $L^2$-projection of exact solution at time $t=t_n$
and $e^n:=\Phi_M^n-\phi^n$ is the difference between the projection $\Phi_M^n$
and the numerical solution $\phi^n$ of the BDF2 implicit scheme \eqref{scheme:PFC BDF2}.
Applying Lemma \ref{lem:Projection-Estimate}, one has
$$\mynorm{\Phi^n-\Phi_M^n}=\mynorm{I_M\brat{\Phi^n-\Phi_M^n}}_{L^2}
\le C_{\phi}\mynorm{I_M\Phi^n-\Phi_M^n}_{L^2}\le C_{\phi}h^{m}\absb{\Phi^n}_{H^{m}}.$$
Once an upper bound of $\mynorm{e^n}$ is available, the claimed error estimate follows immediately,
\begin{align}\label{Triangle-Projection-Estimate}
\mynorm{\Phi^n-\phi^n}\le\mynorm{\Phi^n-\Phi_M^n}+\mynorm{e^n}
\le  C_{\phi}h^{m}\absb{\Phi^n}_{H^{m}}+\mynorm{e^n}
\quad \text{for $1\le n\le N$}.
\end{align}

To bound $\mynorm{e^n}$, we consider two stages:
Stage 1 analyzes the space consistency error for a semi-discrete system having a projected solution $\Phi_M$;
With the help of the DOC kernels $\theta_{k-j}^{(k)}$
and the maximum norm solution estimates in Lemma \ref{lem:Bound-Solution} and \eqref{maximum bound projection},
Stage 2 derives the error estimate from a fully discrete error system
by the standard $L^2$ norm analysis.

\paragraph{Stage 1: Consistency analysis of semi-discrete projection}
A substitution of the projection solution $\Phi_M$ and differentiation  operator $\Delta_h$ into the original equation \eqref{cont: Problem-PFC} yields the semi-discrete system
\begin{align}\label{Projection-Equation}
\partial_t\Phi_M
=\Delta_h\mu_M
+\zeta_{P}\quad\text{with}\quad
\mu_M=(1+\Delta_h)^2\Phi_M + \brab{\Phi_M}^3 -\epsilon\Phi_M,
\end{align}
where $\zeta_{P}(\mathbf{x}_h,t)$ represents the spatial consistency error arising from the projection of exact solution, that is,
\begin{align}\label{Projection-truncation error}
\zeta_{P}:=\partial_t\Phi_M-\partial_t\Phi
+\Delta\mu - \Delta_h\mu_M\quad \text{for $\mathbf{x}_h\in\Omega_{h}$.}
\end{align}
We will bound $\mynorm{\zeta_{P}}$ by applying the triangle inequality,
\begin{align*}
\mynorm{\zeta_{P}}
\le \mynorm{\partial_t\Phi_M-\partial_t\Phi}
+\mynorm{\Delta\mu - \Delta_h\mu_M}
\triangleq\mathbb{I}_1 + \mathbb{I}_2.
\end{align*}
It is easy to check that $I_M\partial_t\Phi_M=\partial_t\Phi_M$ since $\partial_t\Phi_M\in\mathscr{F}_M$,
so one has
\begin{align}
\mathbb{I}_1=&\,\mynorm{\partial_t\bra{\Phi_M-\Phi}}
=\mynormb{I_M\kbra{\partial_t\bra{\Phi_M-\Phi}}}_{L^2}\nonumber\\
\le&\, \mynorm{\partial_t\bra{\Phi_M-\Phi}}_{L^2}
+\mynorm{I_M\partial_t\Phi-\partial_t\Phi}_{L^2}\le C_\phi h^m\mynorm{\partial_t\Phi}_{H^m},\label{Time-Projection}
\end{align}
where Lemma \ref{lem:Projection-Estimate} was used in the second inequality.
It remains to bound the term $\mathbb{I}_2$.
Noticing that $\Delta_h\Phi_M=\Delta\Phi_M$ at the discrete level for $\Phi_M\in\mathscr{F}_{M}$,
we use Lemma \ref{lem:Projection-Estimate} to derive that
\begin{align}
\mynorm{\Delta_h^s\bra{\Phi_M-\Phi}}
&=\mynormb{I_M\kbra{\Delta^s\bra{\Phi_M-\Phi}}}_{L^2}
\le C_\phi h^m\mynorm{\Phi}_{H^{m+2s}}\quad\text{for $s=1,2,3$.} \label{Space-Projection}
\end{align}
For the nonlinear term of $\mathbb{I}_2$,
an application of the triangle inequality leads to
\begin{align}\label{Nonlinear-Projection-1}
\mynorm{\Delta\Phi^3-\Delta_h\Phi_M^3}
\le \mynorm{\Delta\bra{\Phi^3-\Phi_M^3}}
+ \mynorm{\Delta\Phi_M^3-\Delta_h\Phi_M^3}
\triangleq\mathbb{I}_{21} + \mathbb{I}_{22}.
\end{align}
For the term $\mathbb{I}_{21}$, the triangle inequality yields
\begin{align}
\mathbb{I}_{21}
&=\mynorm{\Delta\bra{\Phi^3-\Phi_M^3}}=\mynorm{I_M\bra{\Delta\bra{\Phi^3-\Phi_M^3}}}_{L^2}\nonumber\\
&\le\mynorm{I_M\bra{\Delta\Phi^3}-\Delta\Phi^3}_{L^2}
+\mynorm{\Delta\bra{\Phi^3-\Phi_M^3}}_{L^2}
+\mynorm{\Delta\Phi_M^3-I_M\bra{\Delta\Phi_M^3}}_{L^2}\nonumber\\
&\le C_\phi h^m\mynorm{\Phi^3}_{H^{m+2}}
+ C_\phi h^m\mynorm{\Phi_M^3}_{H^{m+2}}
+\mynorm{\Delta\bra{\Phi^3-\Phi_M^3}}_{L^2}\nonumber\\
&\le  C_\phi h^m\mynorm{\Phi}_{H^{m+2}}^3
+C_\phi h^m\mynorm{\Phi_M}_{H^{m+2}}^3
+\mynorm{\Delta\bra{\Phi^3-\Phi_M^3}}_{L^2},\label{Nonlinear-Projection-2}
\end{align}
in which Lemma \ref{lem:Projection-Estimate}
and the Sobolev embedding inequality have been used in the second and third inequalities,
respectively. For the remainder term in the above inequality \eqref{Nonlinear-Projection-2},
we have the following estimate
\begin{align*}
\mynorm{\Delta\bra{\Phi^3-\Phi_M^3}}_{L^2}&\le
\mynorm{\bra{\Phi^2+\Phi\Phi_M+\Phi_M^2}\bra{\Phi-\Phi_M}}_{H^2}
\le  C_\phi h^m \mynorm{\Phi}_{H^4}^2\mynorm{\Phi}_{H^{m+2}},
\end{align*}
where the estimation \eqref{Projection-Estimate}
and the Sobolev embedding inequality were  used in the second inequality.
Inserting it into  \eqref{Nonlinear-Projection-2}
yields $\mathbb{I}_{21}\le C_\phi h^m$.
In a similar manner, one can find that
$\mathbb{I}_{22}=\mynorm{\Delta\Phi_M^3-\Delta_h\Phi_M^3}\le C_\phi h^m$.
Thus, going back to \eqref{Nonlinear-Projection-1} gives
$\mynorm{\Delta\Phi^3-\Delta_h\Phi_M^3}\le C_\phi h^m$.
The estimations \eqref{Space-Projection}
and \eqref{Nonlinear-Projection-1} yield that
$\mathbb{I}_2=\mynorm{\Delta\mu - \Delta_h\mu_M}\le C_\phi h^m$.

As a consequence,  one has $\mynorm{\zeta_{P}}\le C_\phi h^m$
and $\mynorm{\zeta_{P}(t_j)}\le C_\phi h^m$ for $j\ge1$.
We now apply Lemma \ref{lem:DOC property}(III) to obtain that
\begin{align}\label{Projection-consistency}
\sum_{k=1}^n\mynormb{\Upsilon_{P}^k}\le C_\phi h^m \sum_{k=1}^n\sum_{j=1}^k\theta_{k-j}^{(k)}
\le C_\phi t_nh^m\quad\text{where}\;\; \Upsilon_{P}^k:=\sum_{j=1}^k\theta_{k-j}^{(k)}\zeta_{P}(t_j)\;\;\text{for $k\ge1$.}
\end{align}

\paragraph{Stage 2: $L^2$ norm error of fully discrete system}
From the projection equation \eqref{Projection-Equation},
one can apply the BDF2 formula to obtain the following approximation equation
\begin{align}\label{Discrete-Projection-Equation}
D_2\Phi_M^n
=\Delta_h\mu_M^n +\zeta_{P}^n + \xi_{\Phi}^n\quad\text{with}\quad
\mu_M^n
=(1+\Delta_h)^2\Phi_M^n + \brat{\Phi_M^n}^3 -\epsilon\Phi_M^n,
\end{align}
where $\xi_{\Phi}^n$ is the local consistency error of BDF2 formula,
and $\zeta_{P}^n:=\zeta_{P}(t_n)$ is defined by \eqref{Projection-truncation error}.
Subtracting the full discrete scheme \eqref{scheme:PFC BDF2} from the approximation equation \eqref{Discrete-Projection-Equation},
we have the following error system
\begin{align}\label{Error-Equation}
D_2e^n
&=\Delta_h\kbrab{(1+\Delta_h)^2e^n +f_{\phi}^ne^n}
+\zeta_{P}^n+\xi_{\Phi}^n\quad\text{for $1\le n\le N$,}
\end{align}
where $f_{\phi}^n :=\brat{\Phi_M^n}^2+\Phi_M^n\phi^n+\brat{\phi^n}^2-\epsilon$.
Thanks to the maximum norm solution estimates in Lemma \ref{lem:Bound-Solution} and \eqref{maximum bound projection},
one has
\begin{align}\label{ieq: max bound nonlinear}
\mynormb{f_{\phi}^n}_{\infty}\le c_1^2+c_0c_1+c_0^2+\epsilon.
\end{align}

Multiplying both sides of equation \eqref{Error-Equation} by the DOC kernels $\theta_{k-n}^{(k)}$,
and summing up $n$ from $n=1$ to $k$, we apply the equality
\eqref{orthogonal equality for BDF2} with $v^j=e^j$ to obtain
\begin{align}\label{Error-Equation-DOC}
\diff e^k
&=\sum_{j=1}^k\theta_{k-j}^{(k)}\Delta_h\kbrab{
(1+\Delta_h)^2e^j + f_{\phi}^j e^j}
+\Upsilon_{P}^k+\Xi_{\Phi}^k\quad\text{for $1\le n\le N$,}
\end{align}
where $\Xi_{\Phi}^k$ and $\Upsilon_{P}^k$ are defined by \eqref{BDF2-global consistency} and \eqref{Projection-consistency}, respectively.
Making the inner product of \eqref{Error-Equation-DOC} with $2e^k$,
and summing up the superscript from 1 to $n$, we have the following equality
\begin{align}\label{Error-Equation-Inner}
\mynormb{e^n}^2-\mynormb{e^0}^2
+\sum_{k=1}^n\mynormb{\diff e^k}^2
&=J^n +2\sum_{k=1}^n\myinnerb{\Upsilon_{P}^k+\Xi_{\Phi}^k,e^k}\quad\text{for $1\le n\le N$,}
\end{align}
where $J^n$ is defined by
\begin{align}\label{Error-quadratic forms}
J^n:=&\,2\sum_{k,j}^{n,k}\theta_{k-j}^{(k)}\myinnerb{e^j+2\Delta_he^j+\Delta_h^2e^j+f_{\phi}^je^j,\Delta_h e^k}\nonumber\\
=&\,2\sum_{k,j}^{n,k}\theta_{k-j}^{(k)}\kbra{\myinnerb{f_{\phi}^je^j+2\Delta_he^j,\Delta_h e^k}-\myinnerb{\nabla_he^j,\nabla_he^k}
-\myinnerb{\nabla_h\Delta_he^j,\nabla_h\Delta_h e^k}}.
\end{align}
We are to handle the quadratic form $J^n$.
By applying Lemma \ref{lem:quadr form inequ} with $v^j:=f_{\phi}^je^j$, $w^k:=\Delta_h e^k$ and
$\varepsilon=2\mathcal{M}_r$, one derives that
\begin{align*}
2\sum_{k,j}^{n,k}\theta_{k-j}^{(k)}&\,\myinnerb{f_{\phi}^je^j+2\Delta_he^j,\Delta_h e^k}
=2\sum_{k,j}^{n,k}\theta_{k-j}^{(k)}\myinnerb{f_{\phi}^je^j,\Delta_h e^k}
+4\sum_{k,j}^{n,k}\theta_{k-j}^{(k)}\myinnerb{\Delta_he^j,\Delta_h e^k}\\
\le&\, 4\mathcal{M}_r\sum_{k,j}^{n,k}\theta_{k-j}^{(k)}\myinnerb{f_{\phi}^je^j,f_{\phi}^ke^k}
+5\sum_{k,j}^{n,k}\theta_{k-j}^{(k)}\myinnerb{\Delta_he^j,\Delta_h e^k}\\
\le&\,\sum_{k,j}^{n,k}\theta_{k-j}^{(k)}\kbra{
4\mathcal{M}_r\myinnerb{f_{\phi}^je^j,f_{\phi}^ke^k}+
250\mathcal{M}_r^3\myinnerb{e^j,e^k}
+2\myinnerb{\nabla_h\Delta_he^j,\nabla_h\Delta_h e^k}},
\end{align*}
where the second inequality was obtained by Lemma \ref{lem:inner quad ineq} with $v^j:=e^j$ and
$\varepsilon=2/5$.
Also, Lemma \ref{lem:DOC property} (I) implies that
$-\sum_{k,j}^{n,k}\theta_{k-j}^{(k)}\myinnerb{\nabla_he^j,\nabla_he^k}\le0$.
Then, by applying the Cauchy-Schwarz inequality and
the maximum norm estimate \eqref{ieq: max bound nonlinear}, we obtain from \eqref{Error-quadratic forms} that
\begin{align*}
J^n\le&\,\sum_{k,j}^{n,k}\theta_{k-j}^{(k)}
\kbra{4\mathcal{M}_r\myinnerb{f_{\phi}^je^j,f_{\phi}^ke^k}
+250\mathcal{M}_r^3\myinnerb{e^j,e^k}}
\le c_2\sum_{k,j}^{n,k}\theta_{k-j}^{(k)}\mynormb{e^j}\mynormb{e^k}.
\end{align*}
Therefore, it follows from \eqref{Error-Equation-Inner} that
\begin{align*}
\mynormb{e^n}^2
\le \mynormb{e^0}^2+c_2\sum_{k=1}^n \mynormb{e^k}
\sum_{j=1}^k\theta_{k-j}^{(k)}\mynormb{e^j}
+ 2\sum_{k=1}^n\mynormb{e^k}\mynormb{\Upsilon_{P}^k+\Xi_{\Phi}^k}\quad\text{for $1\le n\le N$.}
\end{align*}

Choosing some integer $n_0$ ($0\le n_0 \le n$) such that
$\mynormb{e^{n_0}}=\max_{0\le k \le n}\mynormb{e^k}$.
Then, taking $n:=n_0$ in the above inequality, one can obtain
\begin{align*}
\mynormb{e^{n_0}}\le \mynormb{e^0}+c_2\sum_{k=1}^{n_0} \mynormb{e^k}
\sum_{j=1}^k\theta_{k-j}^{(k)}+ 2\sum_{k=1}^{n_0}
\mynormb{\Upsilon_{P}^k+\Xi_{\Phi}^k}.
\end{align*}
We know that $\sum_{j=1}^k\theta_{k-j}^{(k)}=\tau_k$ due to Lemma \ref{lem:DOC property}(III).
Thus one gets
\begin{align*}
\mynormb{e^n}\le\mynormb{e^{n_0}}
&\le \mynormb{e^0}+c_2\sum_{k=1}^n \tau_k \mynormb{e^k}+ 2\sum_{k=1}^n\mynormb{\Upsilon_{P}^k+\Xi_{\Phi}^k}.
\end{align*}
Under the time-step size restriction $\tau\le 1/\bra{2c_2}$, we have
\begin{align*}
\mynormb{e^n}\le 2\mynormb{e^0}+2c_2\sum_{k=1}^{n-1} \tau_k \mynormb{e^k}
+ 4\sum_{k=1}^n\mynormb{\Upsilon_{P}^k+\Xi_{\Phi}^k}.
\end{align*}
The discrete Gr\"onwall inequality \cite[Lemma 3.1]{Liao2019Adaptive} yields the following estimate
\begin{align*}
\mynormb{e^n}&\le 2\exp\bra{2c_2t_{n-1}}\braB{\mynormb{e^0}
 + 2\sum_{k=1}^n\mynormb{\Upsilon_{P}^k}+ 2\sum_{k=1}^n\mynormb{\Xi_{\Phi}^k}}\\
 &\le2\exp\bra{2c_2t_{n-1}}\braB{C_\phi h^m+C_\phi t_nh^m
 +2\sum_{k=1}^n\mynormb{\Xi_{\Phi}^k}}\quad\text{for $1\le n\le N$,}
\end{align*}
in which the estimate \eqref{Projection-consistency}
and the initial error $\mynormb{e^0}=\mynormb{\Phi_M^0-\phi^0}\le C_\phi h^m$ have been used.
Moreover,
Lemma \ref{lem:BDF2-Consistency-Error} together with
$\mynormb{\partial_{t}^s\Phi}=\mynormb{I_M\partial_{t}^s\Phi}_{L^2}\le C_\phi\mynormb{\partial_{t}^s\Phi}_{L^2}$ $(s=2,3)$,
due to Lemma \ref{lem:Projection-Estimate}, gives the bound of temporal error term $\sum_{k=1}^n\mynormb{\Xi_{\Phi}^k}$.
Therefore, one obtains the desired estimate from the triangle inequality \eqref{Triangle-Projection-Estimate}
and completes the proof.
\end{proof}

\section{Numerical experiments}
In this section, we apply the variable-step BDF2 scheme \eqref{scheme:PFC BDF2}
to simulate the PFC equation \eqref{cont: Problem-PFC} numerically.
Always, a simple iteration is employed to solve the nonlinear algebra equations
at each time level with the termination error $10^{-12}$.

\subsection{Tests on random time meshes}
\begin{example}\label{PFC-BDF2-Accuracy-Test}
We take $\epsilon=0.02$ and consider the exterior-forced PFC model
$\partial_{t}\Phi=\Delta\mu+ g(\mathbf{x},t)$
in the domain $\Omega=(0,8)^{2}$ such that it has a solution $\Phi=\cos(t)\sin(\frac{\pi}{2}x)\sin(\frac{\pi}{2}y)$.
\end{example}

\begin{table}[htb!]
\begin{center}
\caption{Accuracy of BDF2 scheme \eqref{scheme:PFC BDF2} on random time mesh.}\label{PFC-BDF2-Time-Error} \vspace*{0.3pt}
\def\temptablewidth{0.7\textwidth}
{\rule{\temptablewidth}{0.5pt}}
\begin{tabular*}{\temptablewidth}{@{\extracolsep{\fill}}cccccc}
  $N$   &$\tau$      &$e(N)$     &Order  &$\max r_k$  &$N_1$\\
  \midrule
  20    &8.45e-02	 &1.97e-04	 &--	 &5.84	 &1 \\
  40    &4.80e-02	 &7.73e-05	 &1.65	 &12.22	 &6\\
  80    &2.41e-02	 &1.23e-05	 &2.66	 &746.55 &11\\
  160	&1.27e-02	 &2.94e-06	 &2.24	 &90.35	 &18\\	
  320	&6.51e-03	 &5.97e-07	 &2.39	 &79.85	 &55\\
\end{tabular*}
{\rule{\temptablewidth}{0.5pt}}
\end{center}
\end{table}		

The time accuracy of variable-step BDF2 method \eqref{scheme:PFC BDF2} is examined via random time meshes.
Let the step sizes $\tau_k:=T\sigma_{k}/S$ for $1\leq k\leq N$,
where $\sigma_{k}\in(0,1)$ is the uniformly distributed random number
and $S=\sum_{k=1}^N\sigma_{k}$.
The discrete $L^2$ norm error $e(N):=\|\Phi(T)-\phi^N\|$ is recorded in each run  and the experimental  order of convergence is computed by
$\text{Order}\approx\log\bra{e(N)/e(2N)}/\log\bra{\tau(N)/\tau(2N)}$,
where $\tau(N)$ denotes the maximum time-step size.

The domain $\Omega=(0,8)^2$ is discretized by using $128 \times 128$ mesh
such that the temporal error dominates the spatial error in each run.
We solve  the problem until time $T=1$.
The numerical results are tabulated in Table \ref{PFC-BDF2-Time-Error},
in which we also record the maximum  time-step size $\tau$,
the maximum step ratio and the number (denote by $N_1$ in  Table \ref{PFC-BDF2-Time-Error})
of time levels with the step ratio $r_k\ge (3+\sqrt{17})/2.$
From these data, we observe that the BDF2 scheme is robustly stable
and second-order accuracy on nonuniform time meshes.

\subsection{Numerical comparisons}

\begin{example}
We take the temperature parameter $\epsilon=0.2$ and consider a randomly initial value
$\Phi^0=0.1+0.02\times \mathrm{rand}(\mathbf{x})$ for the PFC model \eqref{cont: Problem-PFC}
in $\Omega=(0,64)^2$, where $\mathrm{rand}(\cdot)$ is the uniformly distributed random number in $(-1,1)$.
The square $\Omega$ is discretized by
a $128\times 128$ uniform mesh.
\end{example}

\begin{figure}[htb!]
\centering
\subfigure[time-step size $\tau=10^{-2}$]{
\includegraphics[width=2.7in,height=1.8in]{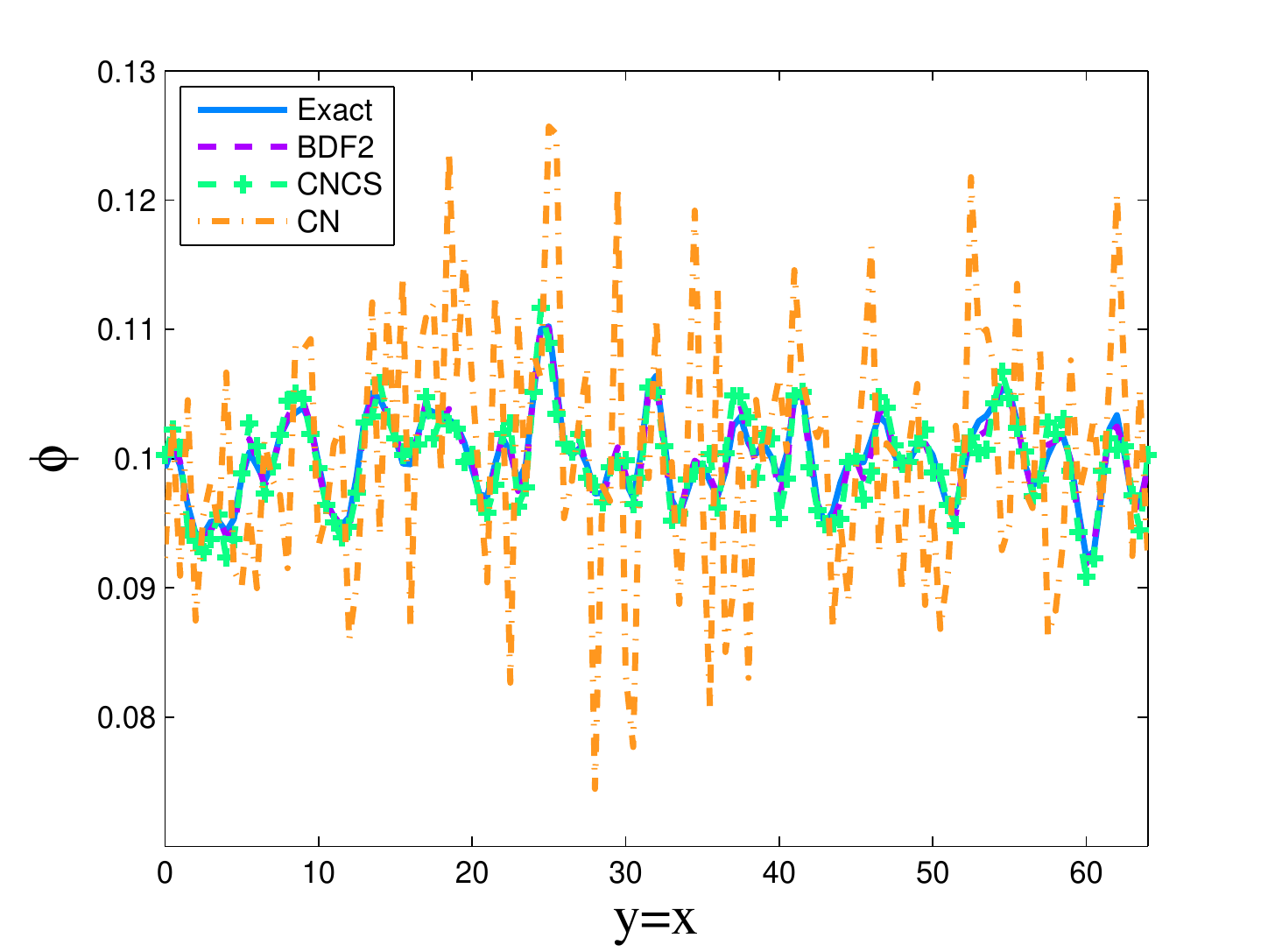}}
\subfigure[time-step size $\tau=10^{-3}$]{
\includegraphics[width=2.7in,height=1.8in]{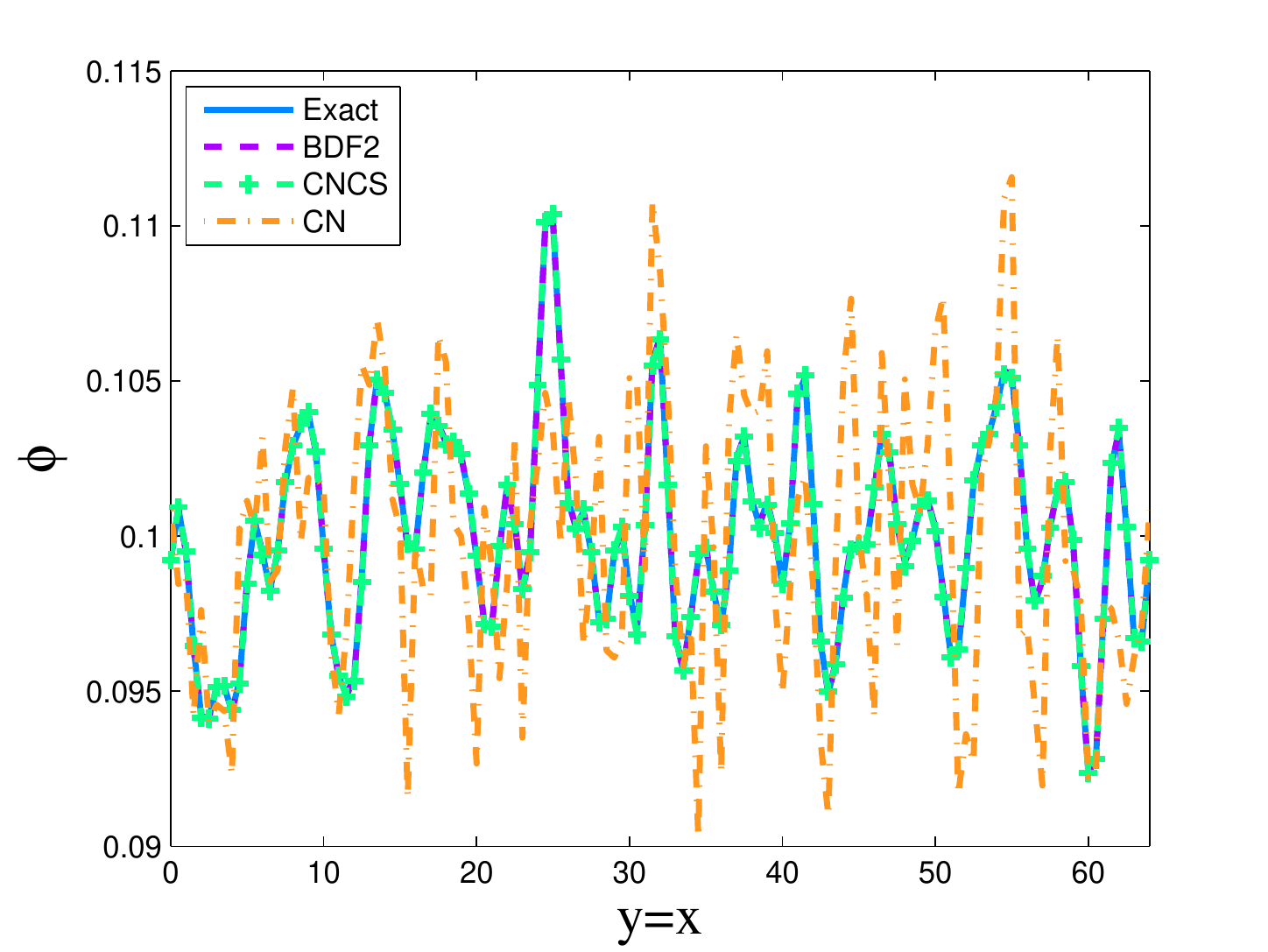}}\\
\subfigure[time-step size $\tau=5\times 10^{-4}$]{
\includegraphics[width=2.7in,height=1.8in]{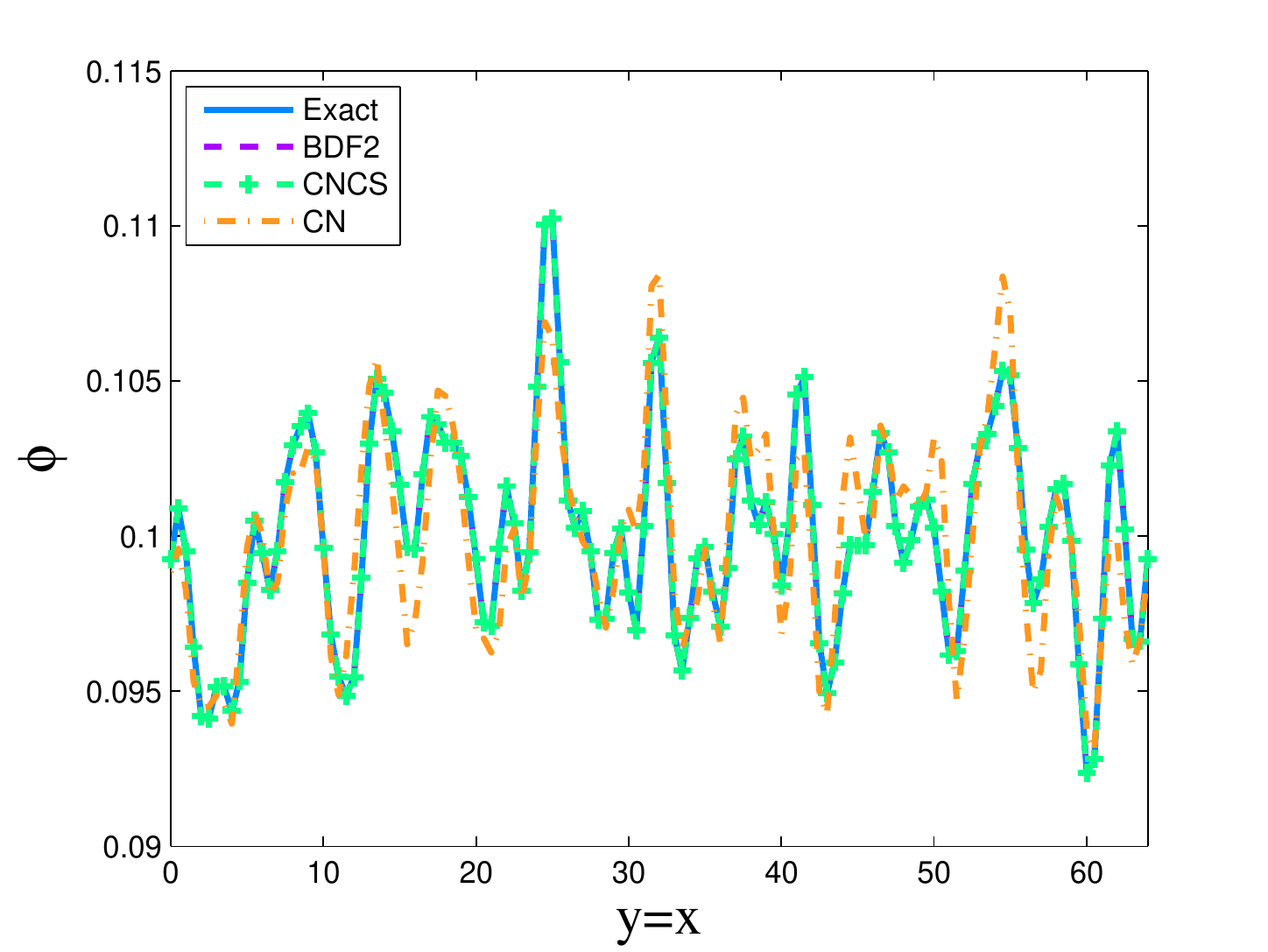}}
\subfigure[time-step size $\tau=2.5\times 10^{-4}$]{
\includegraphics[width=2.7in,height=1.8in]{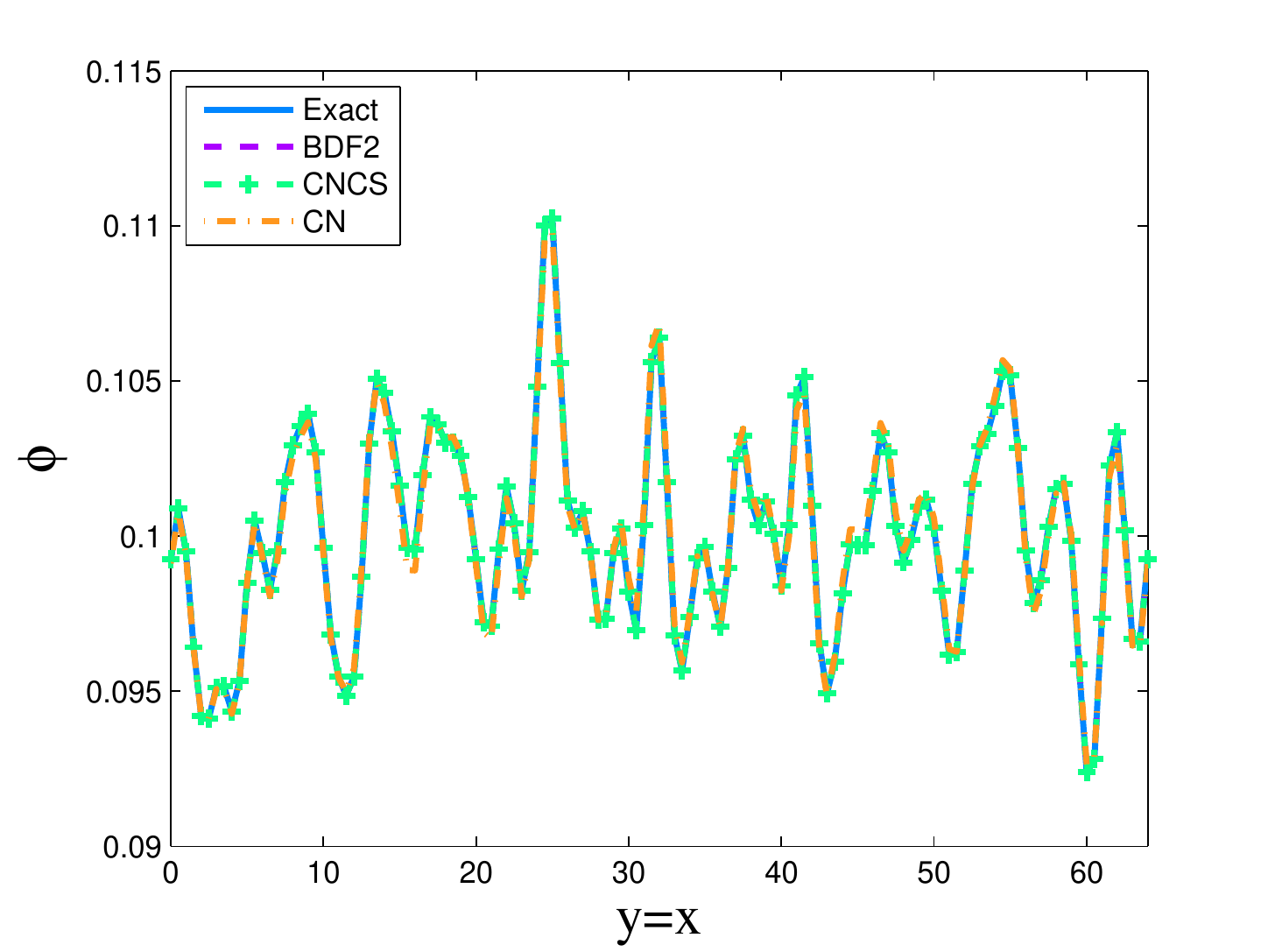}}\\
\caption{Solution curves of BDF2, CN and CNCS methods at the final time $T=0.01$.}
\label{exmaple2:compar BDF2 CN solution}
\end{figure}

To begin with, we examine the numerical behaviors near the initial time
by comparing the BDF2 method \eqref{scheme:PFC BDF2}
with the unconditionally energy stable Crank-Nicoslon (CN) method \cite{zhang2013an},
\begin{align*}
\partial_\tau\phi^n
&=\Delta_h\mu^{n-\frac12},\quad
\mu^{n-\frac12}
=(1+\Delta_h)^2\phi^{n-\frac12}
+\frac{1}{2}\kbrab{(\phi^{n})^2+(\phi^{n-1})^2}\phi^{n-\frac12}
-\epsilon\phi^{n-\frac12},
\end{align*}
and the Crank-Nicoslon convex-splitting (CNCS) scheme  \cite{wise2009an,Dong2018Convergence},
\begin{align*}
\partial_\tau\phi^n
&=\Delta_h\hat{\mu}^{n-\frac12},\quad
\hat{\mu}^{n-\frac12}
=\Delta_h^2\phi^{n-\frac12}
+\Delta_h\hat{\phi}^{n-\frac12}
+\frac{1}{2}\kbrab{(\phi^{n})^2+(\phi^{n-1})^2}\phi^{n-\frac12}
+\bra{1-\epsilon}\phi^{n-\frac12},
\end{align*}
where $\phi^{n-\frac12}:=(\phi^n+\phi^{n-1})/2$ and $\hat{\phi}^{n-\frac12}:=3\phi^{n-1}-\phi^{n-2}$.
We note that the first-order convex-splitting scheme \cite{wise2009an}
is employed to start the CNCS  scheme.
Our computations use a small $T=0.01$ and the reference solution
is computed by the uniform BDF2 method with a vary small time-step size $\tau=10^{-4}$.
The solutions with different time-step sizes are plotted in Figure
\ref{exmaple2:compar BDF2 CN solution}.
In subplot (a), after one step using $\tau=T=10^{-2}$,
the BDF2 solution
is in good with the reference solution,
and the CNCS solution is slightly different from the reference solution,
while the CN solution is  completely different from the reference solution.
Subplot (b) depicts the approximations of 10 steps using $\tau=T/10=10^{-3}$.
We observe that the CN solutions have non-physical oscillations.
The subplots (c)-(d) show the numerical results
after 20 and 40 steps, respectively.
It is seen that the numerical oscillations in the CN solutions are gradually dissipated by very small time-steps.

\begin{table}[htb!]
\begin{center}
\caption{Average iteration numbers and average CPU time (in seconds) at
each time level in BDF2, CN and CNCS methods until $T=5$.}
\label{example2:compar cpu time} \vspace*{0.3pt}
\def\temptablewidth{0.7\textwidth}
{\rule{\temptablewidth}{0.5pt}}
\begin{tabular*}{\temptablewidth}{@{\extracolsep{\fill}}ccccccc}
\multirow{2}{*}{$\tau$} &\multicolumn{2}{c}{BDF2} &\multicolumn{2}{c}{CNCS} &\multicolumn{2}{c}{CN}\\
\cline{2-3}          \cline{4-5}         \cline{6-7}
  &Iter  &CPU  &Iter  &CPU &Iter  &CPU\\
  \midrule
  $10^{-1}$    &5.1224  &0.0122	  &4.3265 &0.0090    &5.1020 &0.0086\\
  $10^{-2}$    &3.9479  &0.0110	  &3.2365 &0.0072	 &4.0100 &0.0072\\
  $10^{-3}$    &3.0064  &0.0095	  &3.0044 &0.0070	 &3.0146 &0.0061\\
\end{tabular*}
{\rule{\temptablewidth}{0.5pt}}
\end{center}
\end{table}	

\begin{figure}[htb!]
\centering
\subfigure[time-step size $\tau=10^{-1}$]{
\includegraphics[width=2.0in]{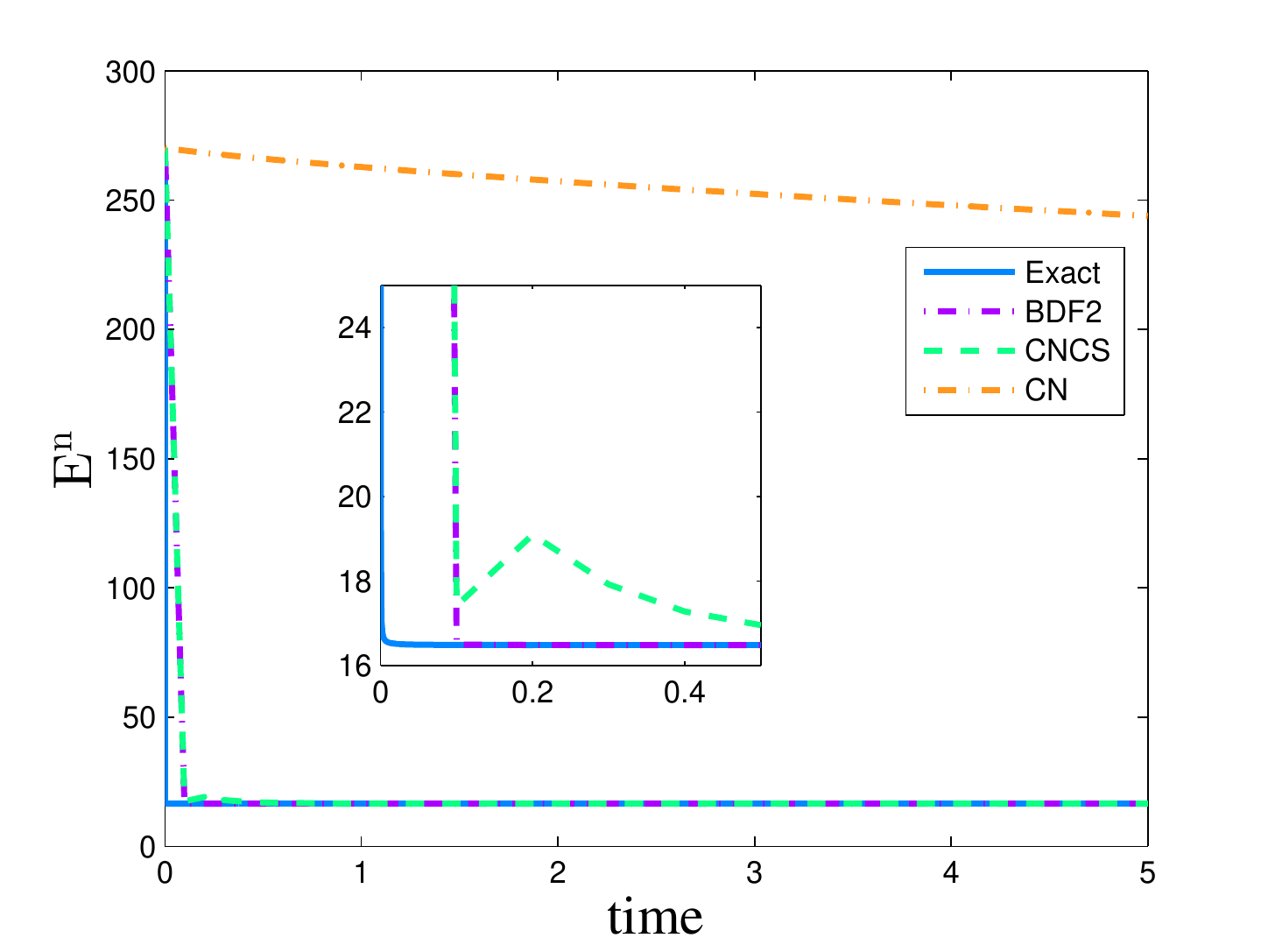}}
\subfigure[time-step size $\tau=10^{-2}$]{
\includegraphics[width=2.0in]{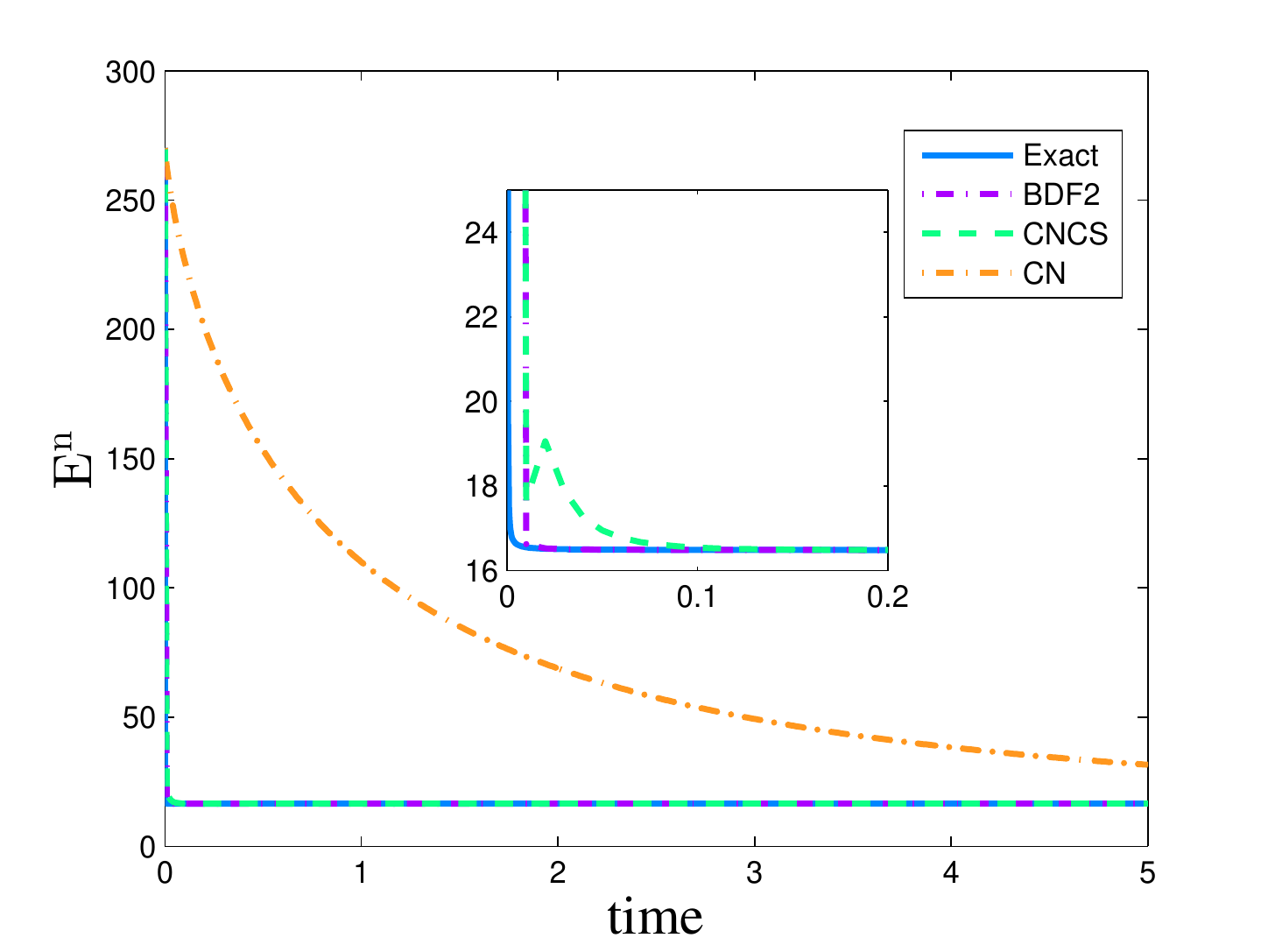}}
\subfigure[time-step size $\tau=10^{-3}$]{
\includegraphics[width=2.0in]{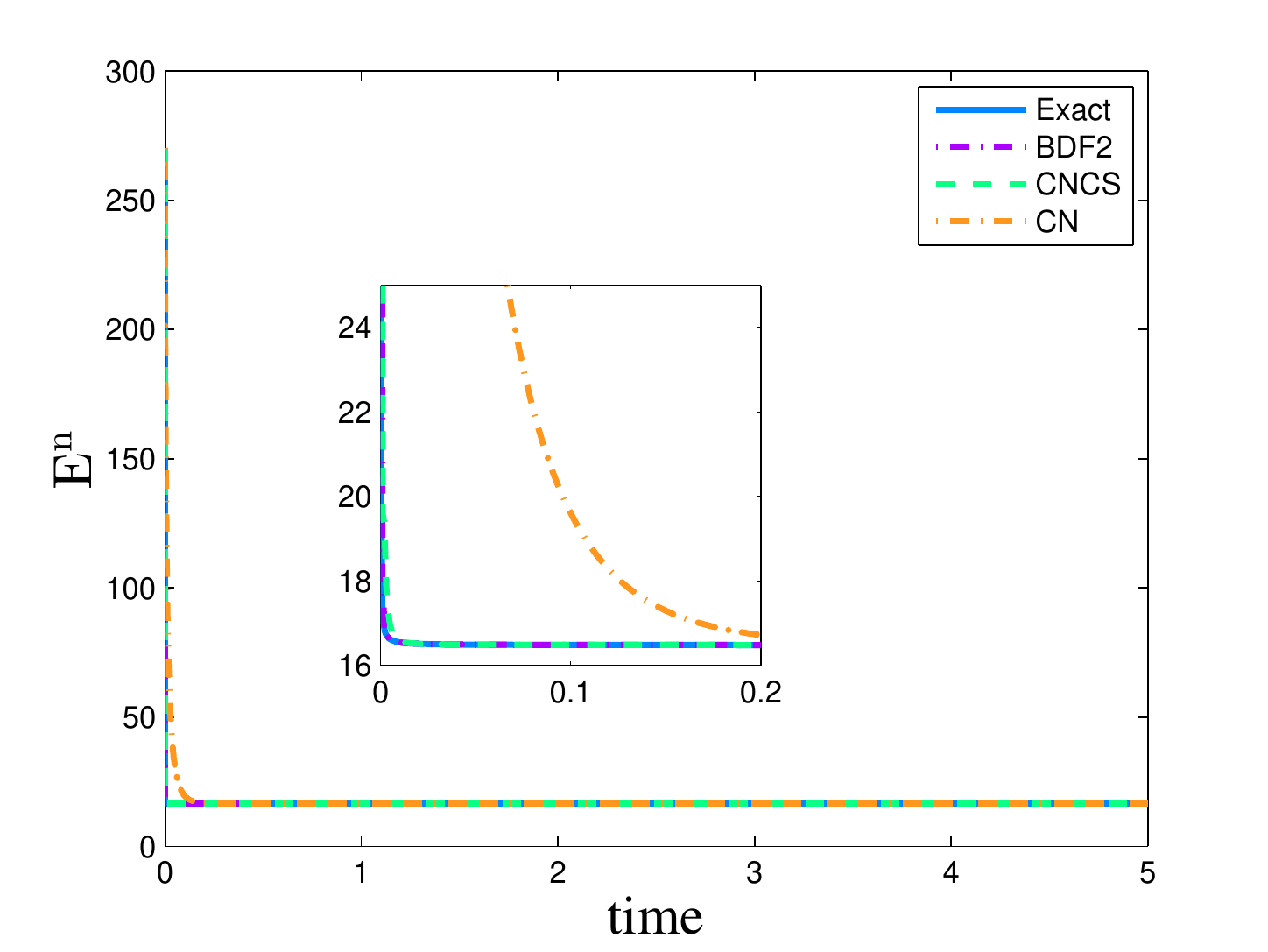}}\\
\caption{Original energy curves of BDF2, CN and CNCS methods until $T=5$.}
\label{example2:compar energy}
\end{figure}

In order to see the numerical performance,
we use the same initial data to
compute the original energy ($E^n=E[\phi^n]$, similarly hereinafter) curves by different time steps until time $T=5$,
see Figure \ref{example2:compar energy}.
The corresponding average iteration numbers (denoted by ``Iter")
and average CPU time (denoted by ``CPU", in seconds) for each time step
are listed in Table \ref{example2:compar cpu time}.
The reference original energy curve is obtained by the uniform BDF2 method with a small time-step $\tau=10^{-4}$.
Table \ref{example2:compar cpu time} shows that the computational cost of BDF2 method is comparable to those
of CN and CNCS  methods. However, as seen in Figure \ref{example2:compar energy},
the original energy curve generated by the CN method deviates from the reference one
and the energy decay property of the CNCS method is numerically destroyed
when some large time-steps are used,
while the BDF2 method generates faithful (original) energy curves for these time-step sizes.

Numerical results indicate that the CN method tends
to generate non-physical oscillations near initial time,
and the BDF2 and CNCS methods can suppress the initial oscillations,
and the former may be a better choice when some large time-step sizes are applied.

\subsection{Adaptive time-stepping strategy}

\begin{algorithm}
\caption{Adaptive time-stepping strategy}
\label{Adaptive-Time-Strategy}
\begin{algorithmic}[1]
\Require{Given $\phi^{n}$ and time step $\tau_{n}$}
\State Compute $\phi^{n+1}$ by using second-order  scheme with time step $\tau_{n}$.
\State Calculate $e_{n+1}=\|\phi^{n+1}-\phi^{n}\|/\|\phi^{n+1}\|$.
\If {$e_{n+1}<tol$ or $\tau_n\le \tau_{\min}$}
 \If {$e_{n+1}<tol$}
\State Update time-step size $\tau_{n+1}\leftarrow\min\{\max\{\tau_{\min},\tau_{ada}\},\tau_{\max}\}$.
 \Else
 \State Update time-step size $\tau_{n+1}\leftarrow \tau_{\min}$.
 \EndIf
\Else
\State Recalculate with time-step size $\tau_{n}\leftarrow\max\{\tau_{\min},\tau_{ada}\}$.
\State Goto 1
\EndIf
\end{algorithmic}
\end{algorithm}

In simulating the phase field problems,
the temporal evolution of phase variables involve
multiple time scales,  such as the growth of a polycrystal discussed
in Example \ref{Growth-Polycrystal}, an initial random perturbation evolves on a fast
time scale, while the later dynamic coarsening evolves on a very slow time scale.
In the following computations, we shall adopt a variant time adaptive  strategy of
 \cite[Algorithm 1]{Gomez2011Provably}
to choose the time step sizes.

The  second-order scheme used in Algorithm \ref{Adaptive-Time-Strategy}
refers to the nonuniform BDF2 scheme in this article.
The adaptive time step $\tau_{ada}$ is given by $\tau_{ada}\bra{e,\tau_{cur}}
=\min\{3.561,\rho\sqrt{{tol}/{e}}\}\tau_{cur},$
where $\rho$ is a default safety coefficient, $tol$ is a reference tolerance,
$e$ is the relative error at each time level,
and $\tau_{cur}$ is the current time step.
In addition,
$\tau_{\max}$ and $\tau_{\min}$ are the predetermined maximum and minimum time steps.
In our computations, if not explicitly specified, we choose the safety coefficient
$\rho=0.9$, the reference tolerance $tol=10^{-3}$,
the maximum time step $\tau_{\max}=0.5$ and the minimum time step $\tau_{\min}=10^{-4}$,
respectively.

\subsection{Growth of a polycrystal}

\begin{example}\label{Growth-Polycrystal}
We take the parameter $\epsilon=0.25$ and use a $256\times 256$ uniform mesh
to discrete the spatial domain $\Omega=(0,256)^2$.
As seeds for nucleation, three random perturbations on the three small square
patches are taken as
$\Phi_0(\mathbf{x})=\bar{\Phi}+A\times \mathrm{rand}(\mathbf{x}),$
where the constant density $\bar{\Phi}=0.285$, $A$ is amplitude and
the random numbers $\mathrm{rand}(\cdot)$ are uniformly distributed in $(-1,1)$.
The centers of three pathes locate at $(128,64),(64,196)$ and $(196,196)$,
with the corresponding amplitudes $A=0.2$, $0.3$ and $0.9$, respectively.
The length of each small square is set to $10$.
\end{example}

\begin{figure}[htb!]
\centering
\includegraphics[width=2.7in,height=1.8in]{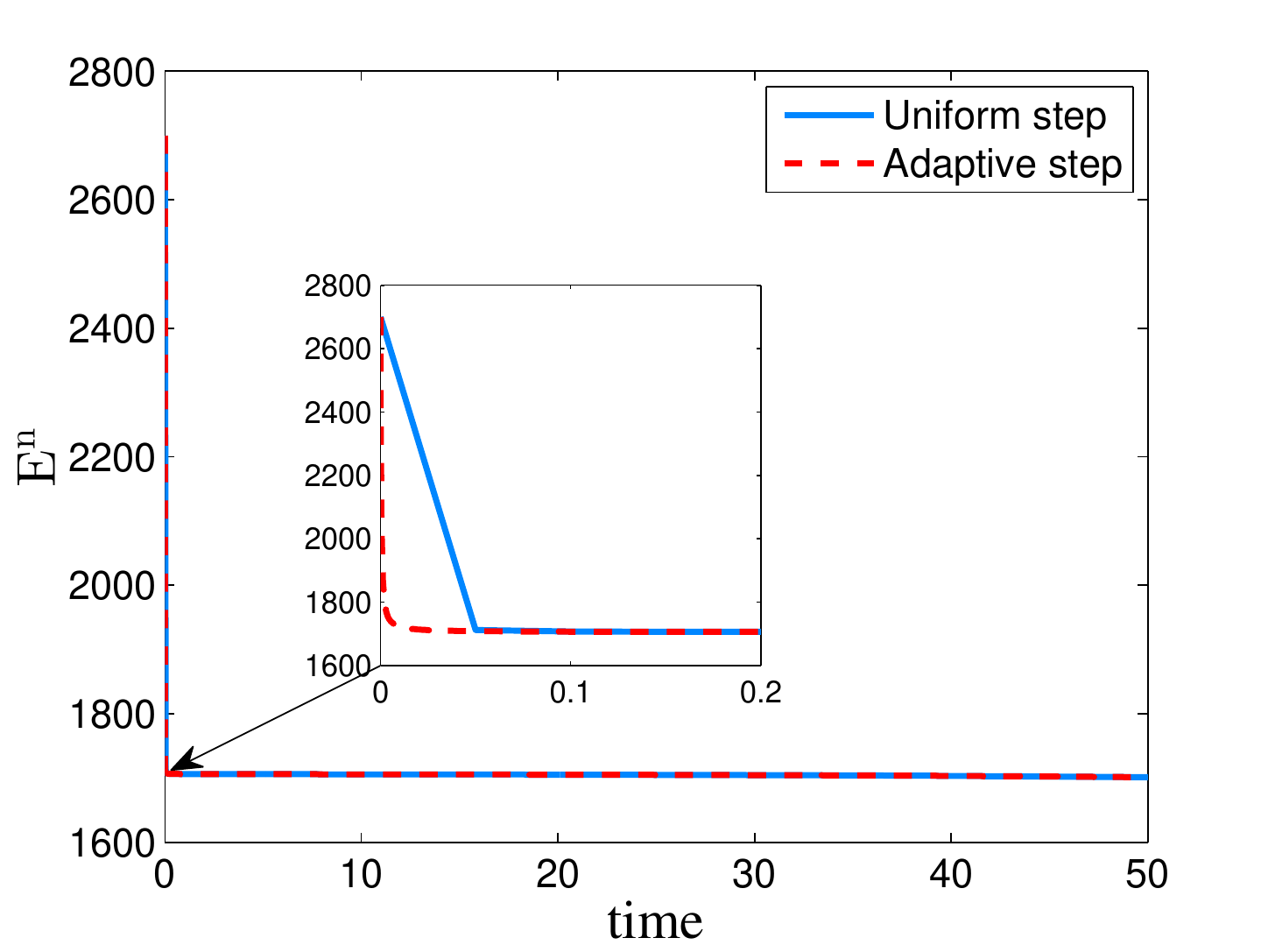}
\includegraphics[width=2.7in,height=1.8in]{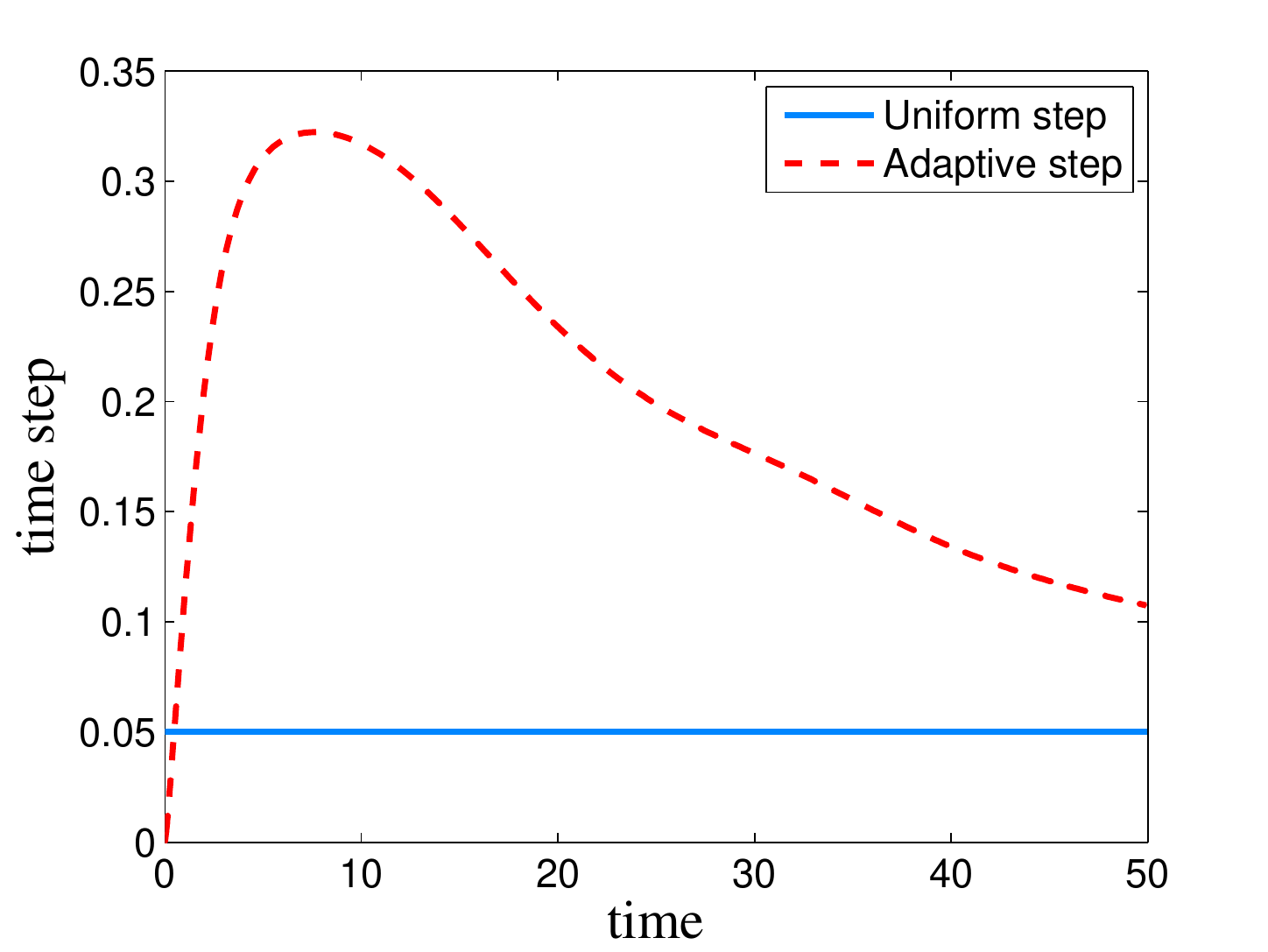}
\caption{Evolutions of original energy (left) and time step sizes (right) of
   the  PFC equation using different time strategies until time $T=50$.}
\label{Comparison-Uniform-Adaptive-Energy}
\end{figure}

\begin{figure}[htb!]
\centering
\includegraphics[width=2.0in]{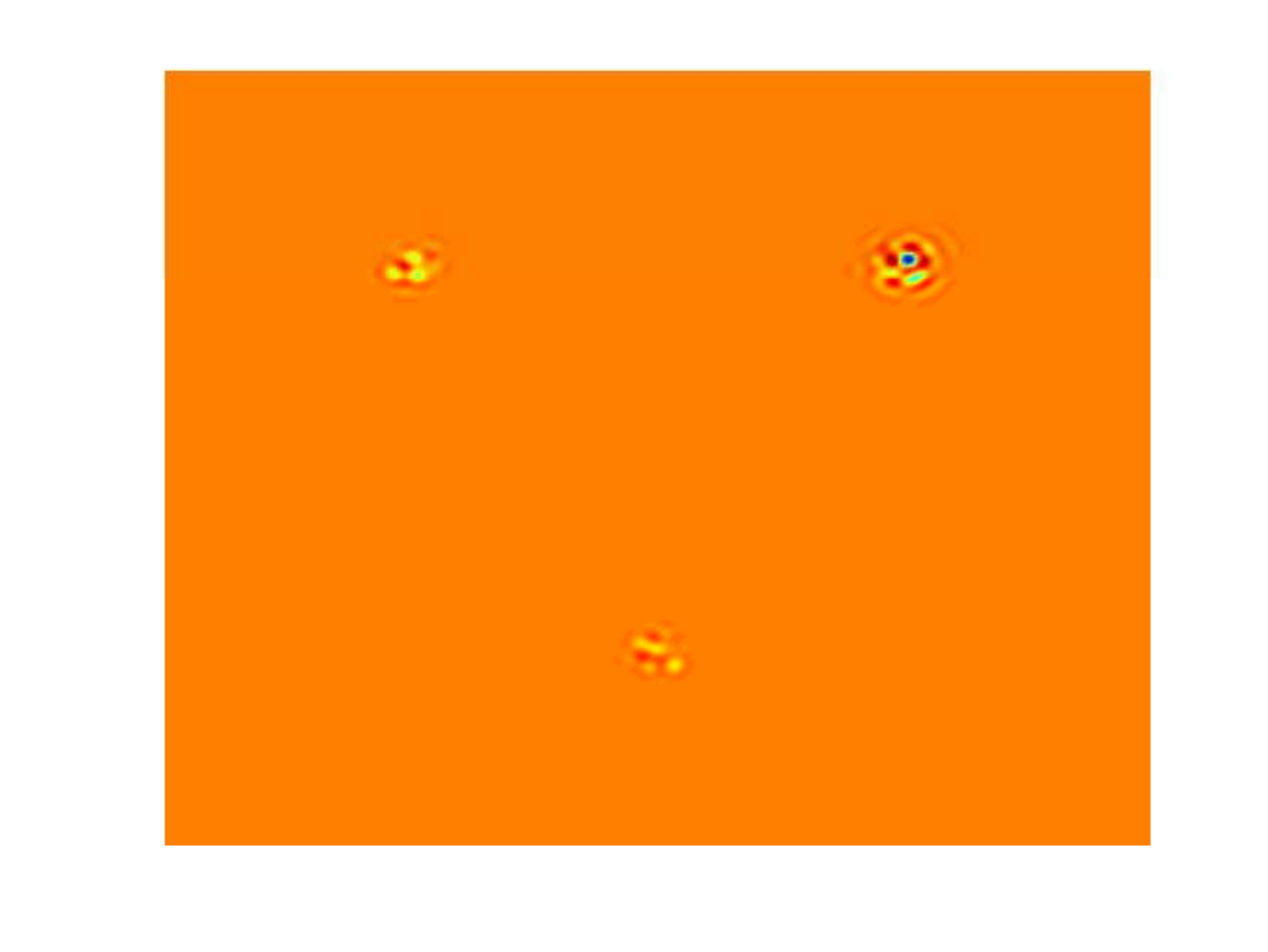}
\includegraphics[width=2.0in]{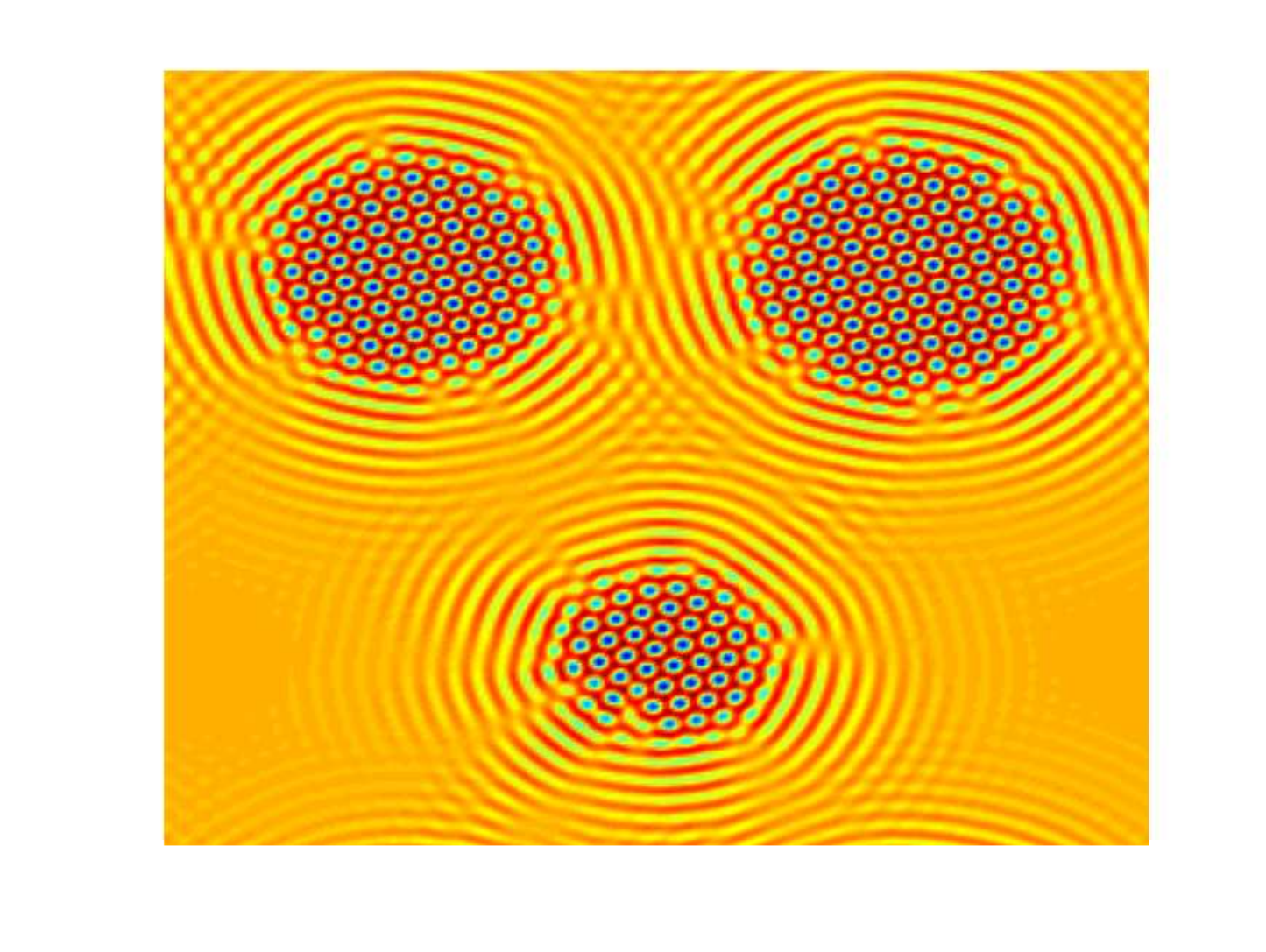}
\includegraphics[width=2.0in]{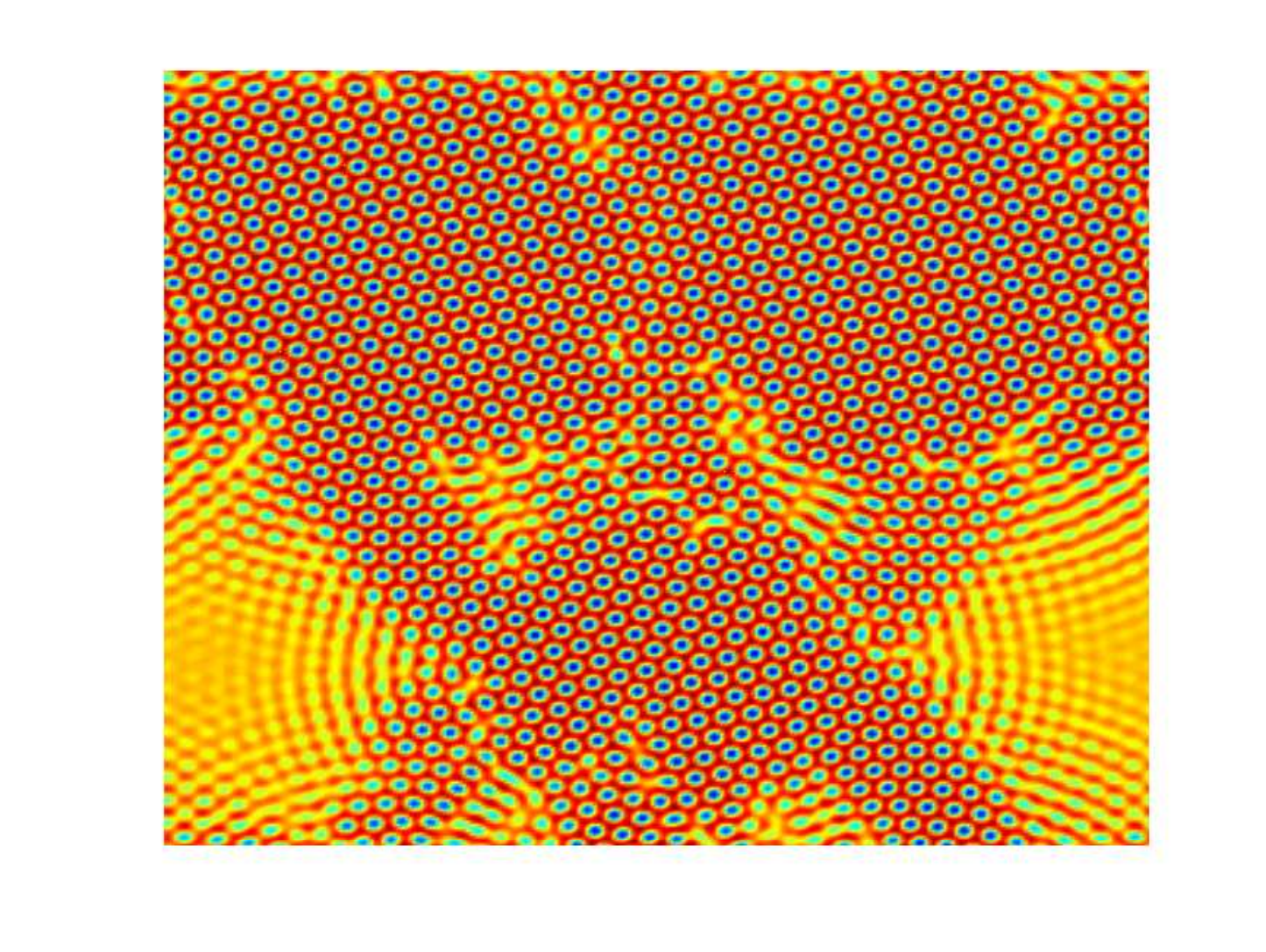}\\
\includegraphics[width=2.0in]{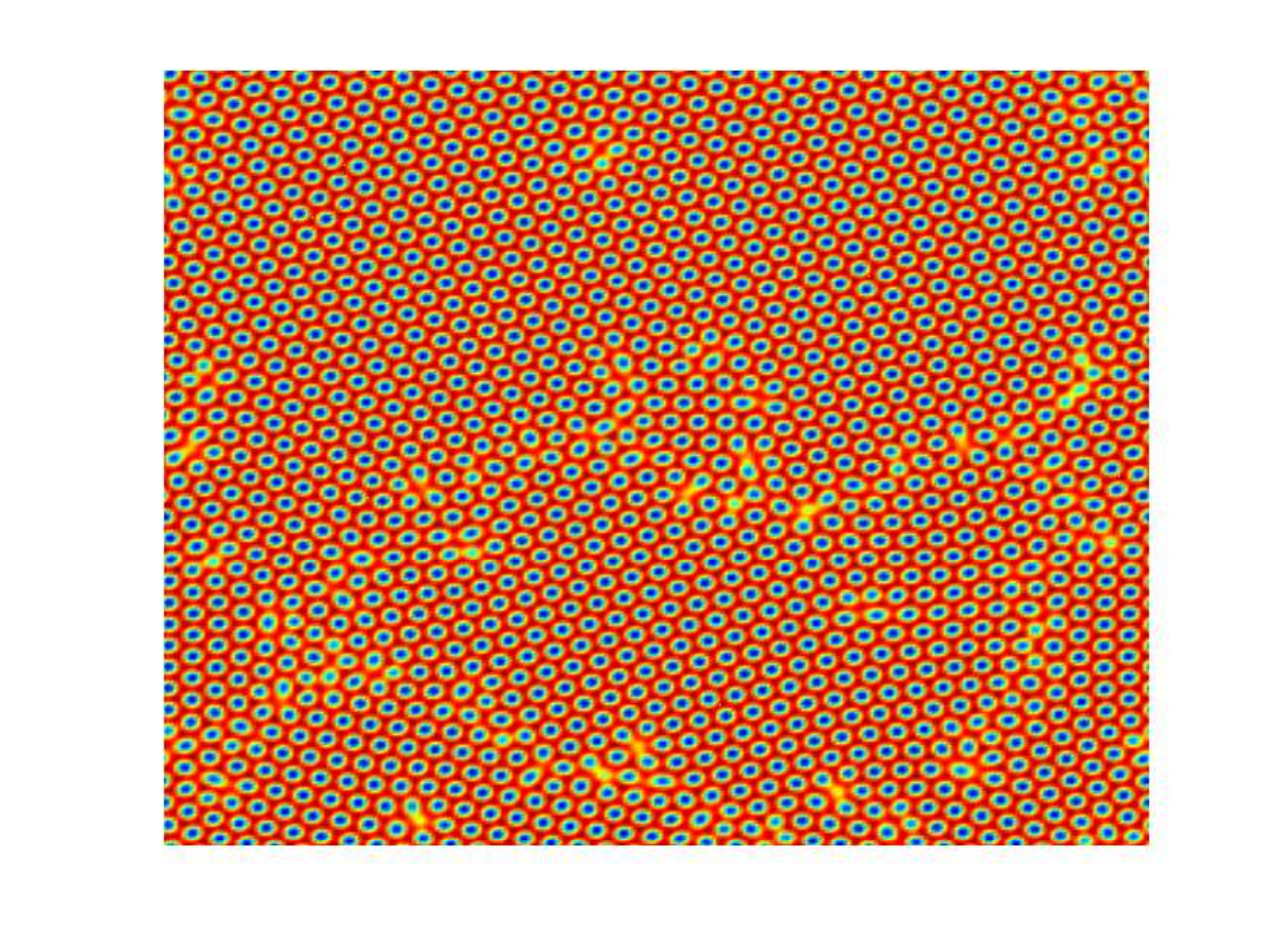}
\includegraphics[width=2.0in]{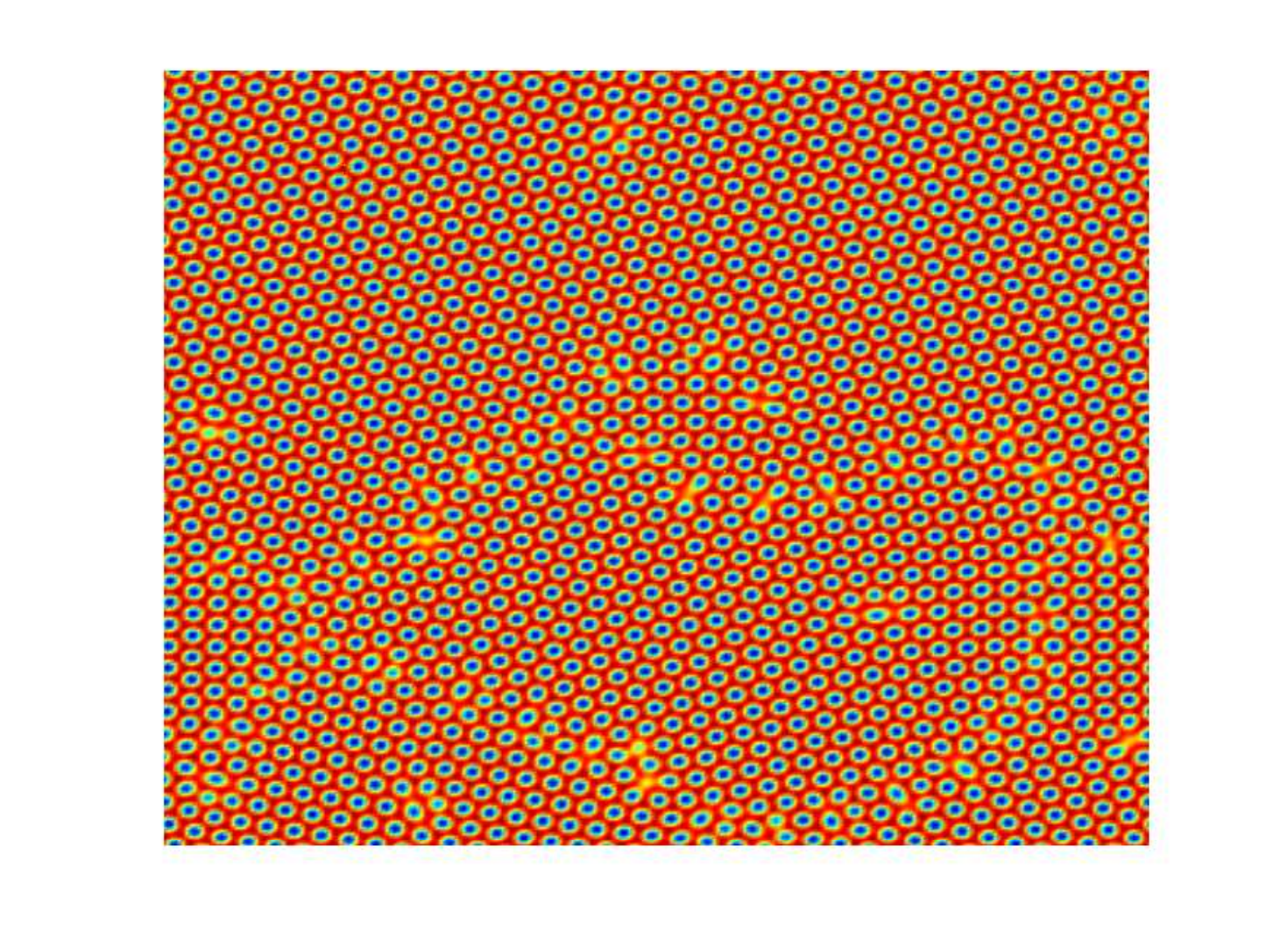}
\includegraphics[width=2.0in]{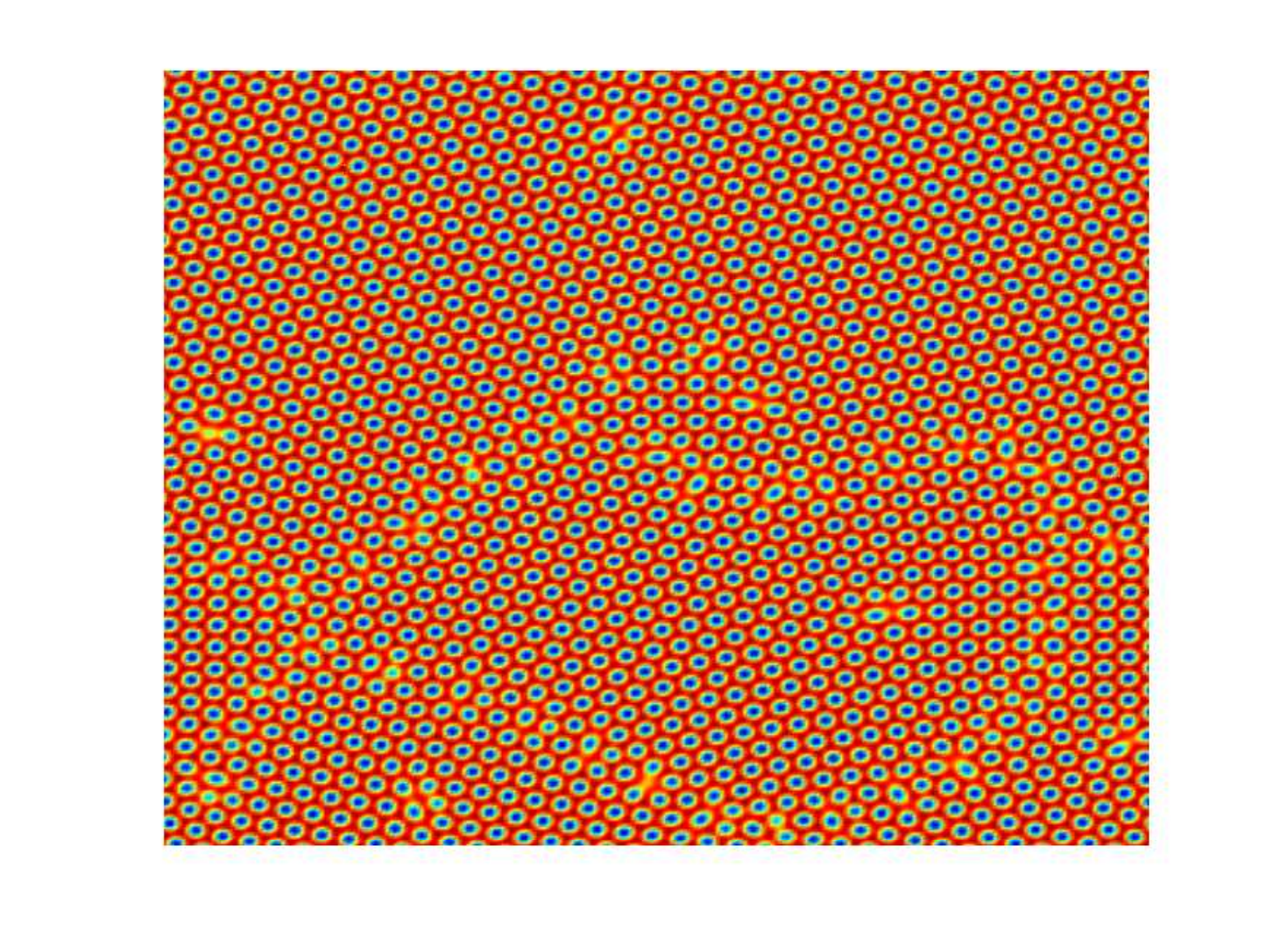}\\
\caption{Solution snapshots of the crystal growth for
   the PFC equation using adaptive time strategy at $t=1, 100, 150, 400, 800,1000$, respectively.}
\label{Snapshots-Dynamics}
\end{figure}

\begin{figure}[htb!]
\centering
\includegraphics[width=2.0in]{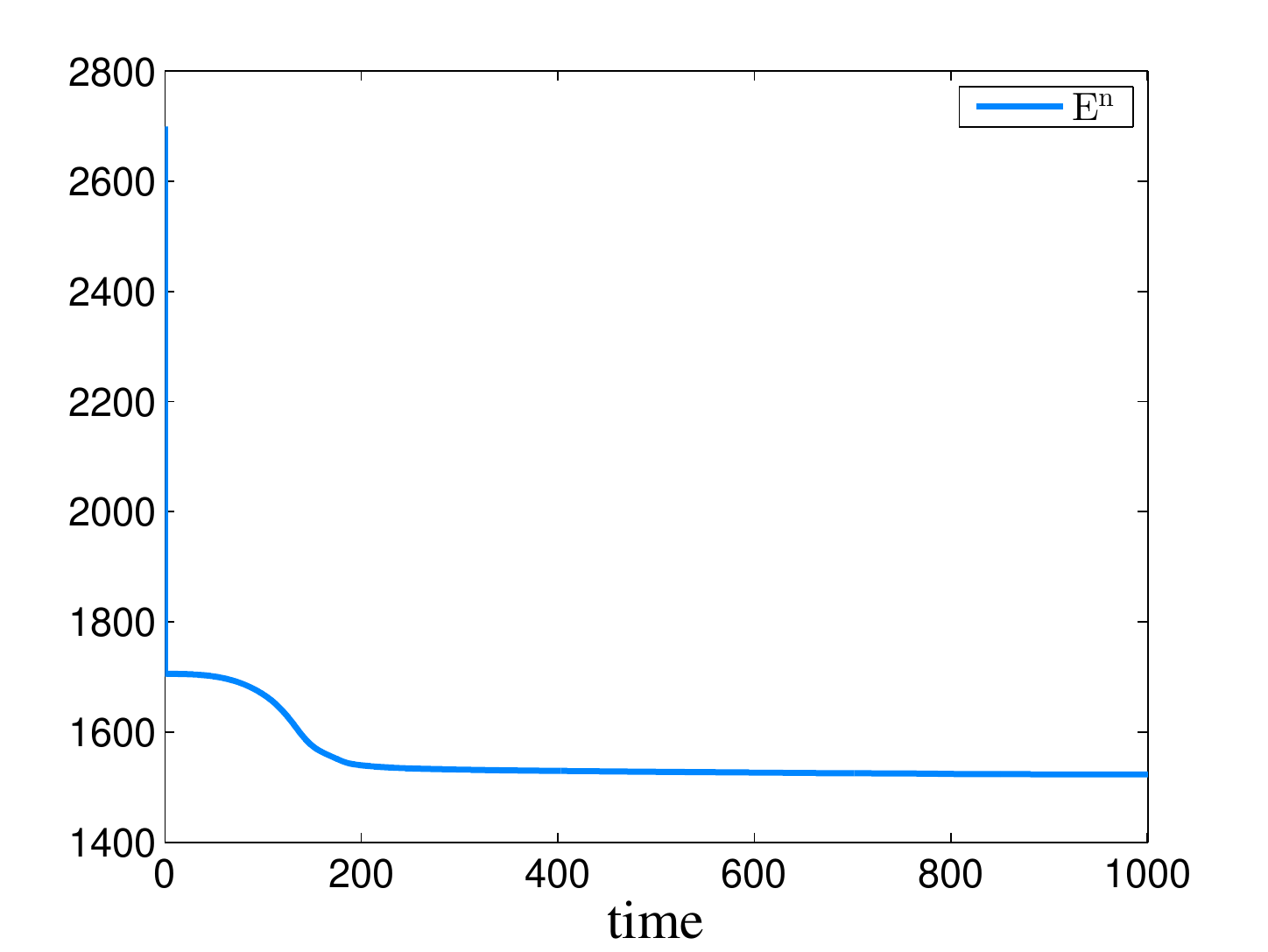}
\includegraphics[width=2.0in]{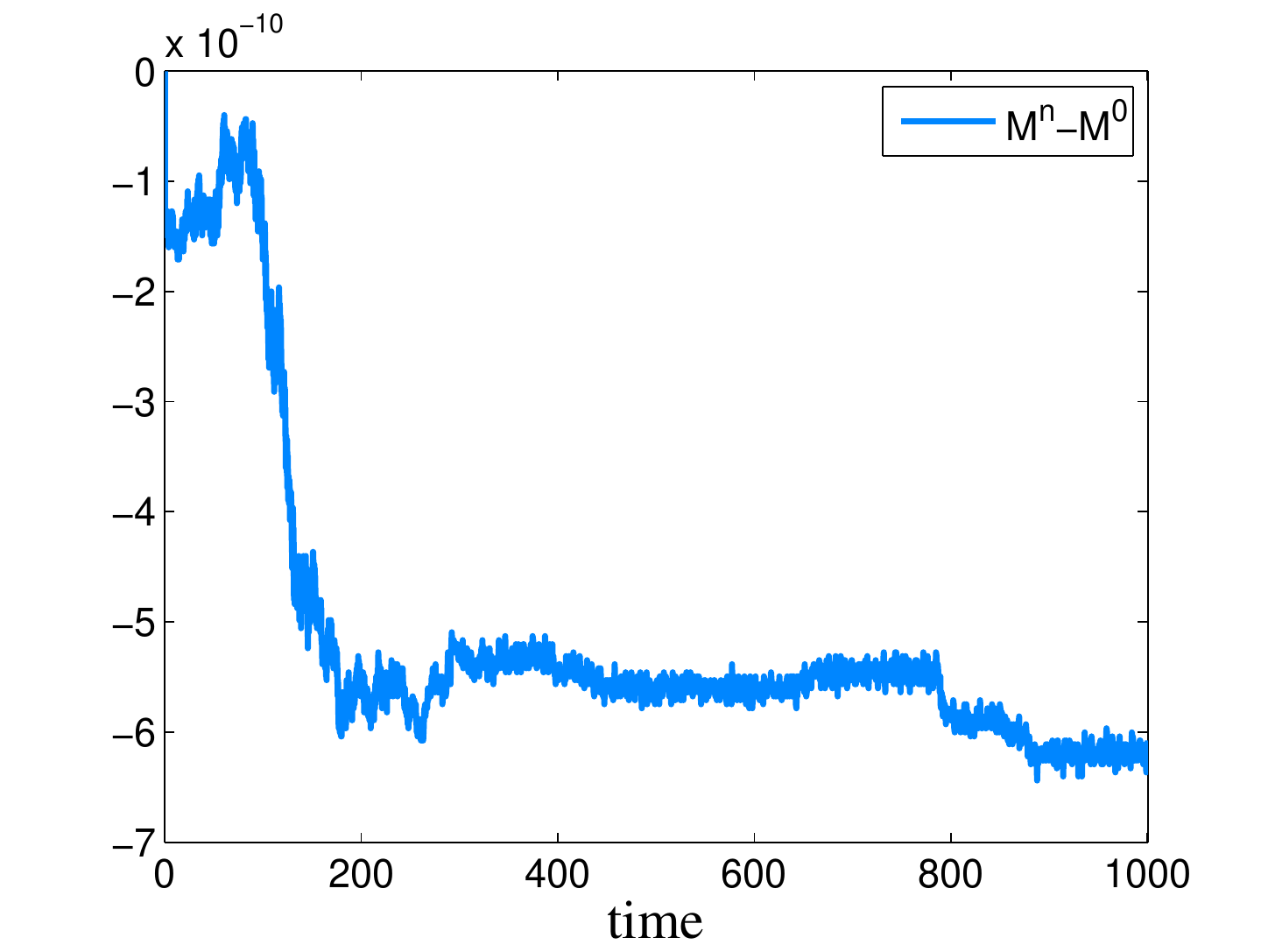}
\includegraphics[width=2.0in]{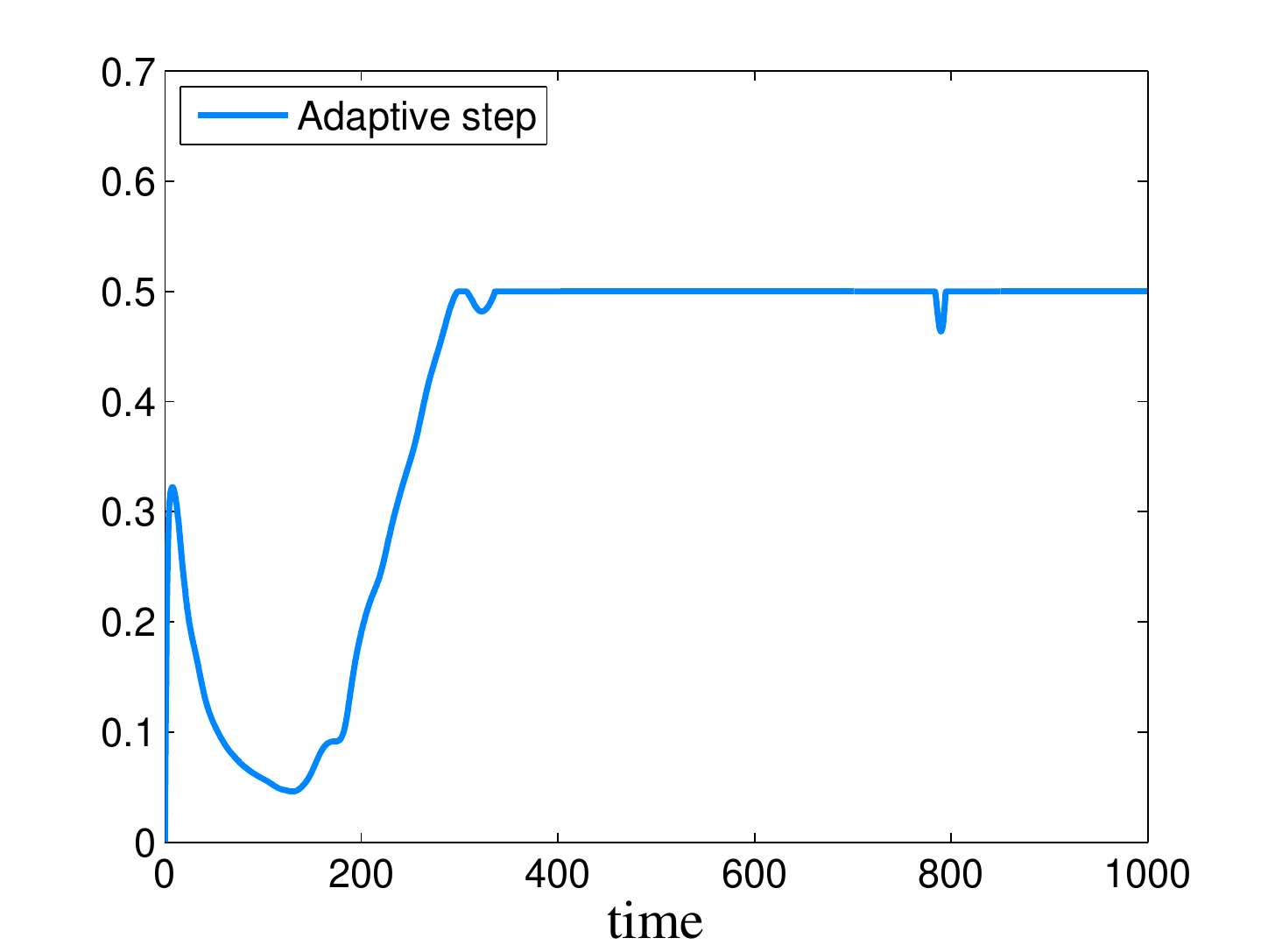}\\
\caption{Evolutions of original energy (left), mass difference (middle) and adaptive
  time steps (right)  of  the crystal growth of PFC equation using adaptive time strategy.}
\label{Adaptive-Energy-Mass}
\end{figure}

We simulate the growth of a polycrystal in a supercooled liquid
with the above random initial liquid density in this example.
We begin with examining the efficiency of adaptive time-stepping  Algorithm \ref{Adaptive-Time-Strategy}
by using different time strategies, i.e., the uniform and adaptive time approaches.
At first, the solution is computed until the time $T=50$ with a constant time step $\tau=0.05$.
We then implement the adaptive strategy described in Algorithm \ref{Adaptive-Time-Strategy}
to simulate the dynamical process by using the same initial data.
The  time evolutions of discrete energies  and the corresponding time-step sizes are plotted
in Figure \ref{Comparison-Uniform-Adaptive-Energy}.
As can be seen, the adaptive energy curve is practically indistinguishable from
that generated by using a small constant step size,
and the former exhibits more details owing to the smaller step sizes are used.
We note that this simulation takes 1000 uniform
time steps with $\tau=0.05$, while the total number of adaptive time steps is 393.
Thus the above numerical results show that the time-stepping adaptive
strategy is computationally efficient.

We now use the BDF2 scheme coupled with Algorithm
\ref{Adaptive-Time-Strategy} to
simulate the growth of a polycrystal in a supercooled liquid.
In the second set of simulations,
we take the time $T=1000$ and the other parameters are the same as given earlier.
The time evolutions of the phase variable are depicted in Figure \ref{Snapshots-Dynamics}.
We see that, the speed of moving interfaces is deeply affected by the initial amplitude,
that is, the larger the amplitude $A$, the faster the polycrystal grows.
Also, three different crystal grains grow and become large enough to form grain boundaries eventually.
The observed phenomena are in good agreement with the published results \cite{li2017an,yang2017linearly}.
The original energy, mass, and adaptive time steps are shown in Figure \ref{Adaptive-Energy-Mass}.
As predicted by our theory, the discrete mass is conservative up to a tolerance of $10^{-10}$.
It is seen that the energy has large variations when the time $t\in [0,200]$,
but it dissipates very slowly when the time escapes.
The right subplot of Figure \ref{Snapshots-Dynamics} shows
that small time step sizes are adopted when the energy dissipates fast, and
large step sizes are utilized when the energy decreases slowly.

\section{Concluding remarks}
\setcounter{equation}{0}
Under the step ratio constraint \textbf{S1},
we proved the  variable-step BDF2 method \eqref{scheme:PFC BDF2} for the PFC model
preserves a discrete modified energy dissipation law,
which implies the maximum norm bound of numerical solution.
The  DOC technique was then improved to establish a concise
convergence analysis of the variable-step BDF2 method. We proved at the first time
that the BDF2 method is convergent
in the $L^2$ norm under the weak step ratio restriction \textbf{S1}.

In our recent work \cite{Liao2020Analysis}, the DOC technique will be further developed
to establish a sharp $L^2$ error estimate on the variable-step BDF2 scheme
for the molecular beam epitaxial growth model without slope selection.
It is expected that the DOC technique would be a useful analysis tool
for other variable-step BDF type methods,
especially when they are combined with the convex splitting technique
or stabilized strategies to achieve unconditionally energy stable in simulating gradient flow problems.
We plan to address these issues in further studies.


\section*{Acknowledgements}
The authors would like to
thank Prof. Xiuling Hu, Prof. Yuezheng Gong,
and Dr. Lin Wang for their valuable discussions and fruitful suggestions.
We also thank  the editor and the anonymous referees for their
valuable comments and suggestions,
which are very helpful for improving the quality of the article.

\appendix

\section{The proof of Lemma \ref{lem:quadr form inequ}}
To facilitate the proof in what follows, we introduce the following matrices
\[
\mathbf{B}_2:=
\left(
\begin{array}{cccc}
b_{0}^{(1)}  &             &            & \\
b_{1}^{(2)} &b_{0}^{(2)}  &            & \\
            &\ddots       &\ddots      &\\
            &             &b_{1}^{(n)} &b_{0}^{(n)}  \\
\end{array}
\right)_{n\times n}\quad\text{and}\quad
\mathbf{\Theta}_2:=
\left(
\begin{array}{cccc}
\theta_{0}^{(1)}  &                  &  & \\
\theta_{1}^{(2)}  &\theta_{0}^{(2)}  &  & \\
\vdots           &\vdots           &\ddots  &\\
\theta_{n-1}^{(n)}&\theta_{n-2}^{(n)}&\cdots  &\theta_{0}^{(n)}  \\
\end{array}
\right)_{n\times n},
\]
where the discrete kernels
$b_{n-k}^{(n)}$ and $\theta_{n-k}^{(n)}$ are defined by \eqref{def:BDF2-kernels} and \eqref{DOC-Kernels}, respectively.
It follows from the discrete orthogonal
identity  \eqref{orthogonal identity} that
\begin{align}\label{matrix: orthogonal identity}
\mathbf{\Theta}_2=\mathbf{B}_2^{-1}.
\end{align}
If the step ratios condition \textbf{S1} holds, Lemma \ref{lem:conv kernels positive} shows that
the real symmetric matrix
\begin{align}\label{matrix: B}
\mathbf{B}:=\mathbf{B}_2+\mathbf{B}_2^T\quad\text{is positive definite,}
\end{align}
that is,
$\myvec{w}^T\mathbf{B}\myvec{w}=2\sum_{k=1}^n w_k \sum_{j=1}^k b_{k-j}^{(k)}w_j\ge
\sum_{k=1}^nR(r_k,r_{k+1})w_k^2/\tau_k$,
where $R(z,s)$ is defined by \eqref{def: step-ratios function} and $\myvec{w}:=\bra{w_1,w_2,\cdots,w_n}^T$.
According to Lemma \ref{lem:DOC property} (I), the real symmetric matrix
\begin{align}\label{matrix: Theta}
\mathbf{\Theta}:=\mathbf{\Theta}_2+\mathbf{\Theta}_2^T\quad\text{is positive definite,}
\end{align}
in the sense of $\myvec{w}^T\mathbf{\Theta}\myvec{w}=2\sum_{k=1}^n w_k \sum_{j=1}^k \theta_{k-j}^{(k)}w_j>0$.

Moreover, we define a diagonal matrix
$\Lambda_\tau:=\text{diag}\bra{\sqrt{\tau_1},\sqrt{\tau_2},\cdots,\sqrt{\tau_n}}$
and
\begin{align}\label{matrix: modified-B2}
\widetilde{\mathbf{B}}_2:=\Lambda_\tau\mathbf{B}_2\Lambda_\tau=
\left(
\begin{array}{cccc}
\tilde{b}_{0}^{(1)}  &             &            & \\
\tilde{b}_{1}^{(2)} &\tilde{b}_{0}^{(2)}  &            & \\
            &\ddots       &\ddots      &\\
            &             &\tilde{b}_{1}^{(n)} &\tilde{b}_{0}^{(n)}  \\
\end{array}
\right)_{n\times n},
\end{align}
where the discrete kernels $\tilde{b}_{0}^{(k)}$ and $\tilde{b}_{1}^{(k)}$ are given by
($r_1\equiv0$)
\[
\tilde{b}_{0}^{(k)}=\frac{1+2r_k}{1+r_k}\quad
\text{and}\quad
\tilde{b}_{1}^{(k)}=-\frac{r_k^{3/2}}{1+r_k}
\quad\text{for $1\leq k\leq n$}.
\]
Some results on the matrix $\widetilde{\mathbf{B}}_2$ are presented as follows.

\begin{lemma}\label{lem: tilde B-positiveDefinite}
If the step ratios condition \textbf{S1} holds, then the minimum eigenvalue of the real symmetric matrix
\begin{align}\label{matrix: modified-B}
\widetilde{\mathbf{B}}:=\widetilde{\mathbf{B}}_2+\widetilde{\mathbf{B}}_2^T
\end{align}
can be bounded by
\[
\lambda_{\min}\brab{\widetilde{\mathbf{B}}}
\ge\min_{1\le k\le n}R_L\bra{r_{k},r_{k+1}}\ge21/40,
\]
where $R_L\bra{z,s}$ is defined by
\begin{align}\label{def: step-ratios function RL}
R_L\bra{z,s}:=\frac{2+4z-z^{3/2}}{1+z}-\frac{s^{3/2}}{1+s}
\quad\text{for $0\le z,s< r_{\mathrm{sup}}=\tfrac{3+\sqrt{17}}{2}.$}
\end{align}
Thus $\widetilde{\mathbf{B}}$ is positive definite and there exists a non-singular upper triangular matrix $\mathbf{L}$ such that
$$\widetilde{\mathbf{B}}=\Lambda_\tau\mathbf{B}\Lambda_\tau=\mathbf{L}^T\mathbf{L}\quad\text{or}\quad \mathbf{B}=\bra{\mathbf{L}\Lambda_\tau^{-1}}^T\mathbf{L}\Lambda_\tau^{-1}.$$
\end{lemma}
\begin{proof}
Note that, $\partial R_L/\partial z=\frac{(1-\sqrt{z})(z+\sqrt{z}+4)}{2(1+z)^2}$.
Thus $R_L\bra{z,s}$ is increasing in $(0,1)$
and decreasing in $(1, r_{\mathrm{sup}})$ with respect to $z$.
Also, $R_L\bra{z,s}$ is decreasing with respect to $s$ such
that $R_L\bra{z,s}< R_L\bra{z,0}$ for any $s\in(0,r_{\mathrm{sup}})$.
Simple calculations show that
\begin{align}\label{lemproof: tilde B-1}
R_L\bra{z,s}\ge\min\big\{R_L\bra{0,r_{\mathrm{sup}}},R_L\bra{r_{\mathrm{sup}},r_{\mathrm{sup}}}\big\}
>21/40\quad\text{for $0\le z,s< r_{\mathrm{sup}}.$}
\end{align}
For any fixed index $n$, by using the definition \eqref{matrix: modified-B2} of $\widetilde{\mathbf{B}}_2$ and the well--known
Gerschgorin's circle theorem, we find  that the minimum eigenvalue of $\widetilde{\mathbf{B}}$ can be bounded by
\begin{align*}
\lambda_{\min}\brab{\widetilde{\mathbf{B}}}
\ge&\,\min_{1\le k\le n-1}
\Big\{2\tilde{b}_{0}^{(k)}-\absb{\tilde{b}_{1}^{(k)}}-\absb{\tilde{b}_{1}^{(k+1)}},
2\tilde{b}_{0}^{(n)}-\absb{\tilde{b}_{1}^{(n)}}\Big\}\\
=&\,\min_{1\le k\le n-1}
\Big\{2\tilde{b}_{0}^{(k)}+\tilde{b}_{1}^{(k)}+\tilde{b}_{1}^{(k+1)},
2\tilde{b}_{0}^{(n)}+\tilde{b}_{1}^{(n)}\Big\}\\
=&\,\min_{1\le k\le n-1}\Big\{R_L\bra{r_{k},r_{k+1}},R_L\bra{r_n,0}\Big\}
\ge \min_{1\le k\le n}R_L\bra{r_{k},r_{k+1}}>21/40,
\end{align*}
where the last estimate follows from \eqref{lemproof: tilde B-1}. It also says that
the real symmetric matrix $\widetilde{\mathbf{B}}$ is positive definite.
Then we complete the proof by noticing the definition \eqref{matrix: B} and
applying the standard Cholesky decomposition of $\widetilde{\mathbf{B}}$.
\end{proof}

\begin{lemma}\label{lem: tilde B2TB2-positiveDefinite}
If the step ratios condition \textbf{S1} holds, the maximum eigenvalue of the real symmetric matrix $\widetilde{\mathbf{B}}_2^T\widetilde{\mathbf{B}}_2$
can be bounded by
\[
\lambda_{\max}\brab{\widetilde{\mathbf{B}}_2^T\widetilde{\mathbf{B}}_2}
\leq\max_{1\le k\le n}R_U\bra{r_{k},r_{k+1}}<R_U\bra{r_{\mathrm{sup}},r_{\mathrm{sup}}}<53/5,
\]
where $\widetilde{\mathbf{B}}_2$ is defined in \eqref{matrix: modified-B2} and $R_U\bra{z,s}$ is defined by
\begin{align}\label{def: step-ratios function RU}
R_U\bra{z,s}:=\frac{(1+2z)(1+2z+z^{3/2})}{\bra{1+z}^2}
+\frac{s^{3/2}(1+2s+s^{3/2})}{\bra{1+s}^2}
\quad\text{for $0\le z,s< r_{\mathrm{sup}}.$}
\end{align}
\end{lemma}
\begin{proof}Obviously, $R_U\bra{z,s}$ is increasing with respect to the two variables $z$ and $s$.
We have $R_U\bra{z,s}<R_U\bra{r_{\mathrm{sup}},r_{\mathrm{sup}}}<53/5$.
From the definition \eqref{matrix: modified-B2} of $\widetilde{\mathbf{B}}_2$, one has
\[
\widetilde{\mathbf{B}}_2^T\widetilde{\mathbf{B}}_2=
\left(
\begin{array}{ccccc}
d_{0}^{(1)} &d_{1}^{(2)} &              &              &\\
d_{1}^{(2)} &d_{0}^{(2)} &d_{1}^{(3)}   &              & \\
            &\ddots      &\ddots        &\ddots        &\\
            &            &d_{1}^{(n-1)} &d_{0}^{(n-1)} &d_{1}^{(n)}\\
            &            &              &d_{1}^{(n)}   &d_{0}^{(n)}  \\
\end{array}
\right)_{n\times n},
\]
where the discrete kernels $d_{0}^{(k)}$ and $d_{1}^{(k)}$ are given by
($r_1\equiv0$)
\[
d_{0}^{(k)}=\braB{\frac{1+2r_k}{1+r_k}}^2
+\frac{r_{k+1}^{3}}{(1+r_{k+1})^2}
\quad\text{and}\quad
d_{1}^{(k)}=-\frac{r_k^{3/2}\bra{1+2r_k}}{\bra{1+r_k}^2}
\quad\text{for $1\leq k\leq n$}.
\]
For any fixed index $n$, the Gerschgorin circle theorem
gives an upper bound of the maximum eigenvalue of
the real symmetric matrix $\widetilde{\mathbf{B}}_2^T\widetilde{\mathbf{B}}_2$, that is,
\begin{align*}
\lambda_{\max}\brab{\widetilde{\mathbf{B}}_2^T\widetilde{\mathbf{B}}_2}
\le&\,\max_{1\le k\le n-1}
\Big\{d_{0}^{(k)}-d_{1}^{(k)}-d_{1}^{(k+1)},
d_{0}^{(n)}-d_{1}^{(n)}\Big\}\\
=&\,\max_{1\le k\le n-1}\Big\{R_U\bra{r_{k},r_{k+1}},R_U\bra{r_n,0}\Big\}
\le \max_{1\le k\le n}R_U\bra{r_{k},r_{k+1}}.
\end{align*}
It completes the proof.
\end{proof}

\begin{lemma}\label{lem:quadr form inequ B}
If \textbf{S1} holds, then the positive definite matrix $\mathbf{\Theta}=(\mathbf{B}_2^{-1})^T\mathbf{B}\mathbf{B}_2^{-1}$ and
\begin{align*}
\sum_{k=1}^n  \sum_{j=1}^k\theta_{k-j}^{(k)} w_k v_j
\le
\frac{\varepsilon}{2}\myvec{v}^T\mathbf{\Theta}\myvec{v}
+\frac{1}{2\varepsilon}\myvec{w}^T\mathbf{B}^{-1}\myvec{w}
\quad \text{for $\varepsilon > 0$}
\end{align*}
for any real vectors $\myvec{v}:=\bra{v_1,v_2,\cdots,v_n}^T$
and $\myvec{w}:=\bra{w_1,w_2,\cdots,w_n}^T$.
\end{lemma}

\begin{proof}
For any fixed index $n$, let $\myvec{u}:=\mathbf{\Theta}_2\myvec{v}$.
The equality \eqref{matrix: orthogonal identity} gives $\myvec{v}=\mathbf{B}_2\myvec{u}$.
In the element-wise sense, one has $u_k=\sum_{j=1}^k\theta_{k-j}^{(k)} v_j$ and $v_k=\sum_{j=1}^kb_{k-j}^{(k)} u_j$.
Then we have
\begin{align}\label{lemp:conv identity}
\myvec{v}^T\mathbf{\Theta}\myvec{v}=2\sum_{k=1}^n  \sum_{j=1}^k\theta_{k-j}^{(k)} v_k v_j
=2\sum_{k=1}^n  \sum_{j=1}^kb_{k-j}^{(k)} u_k u_j=\myvec{u}^T\mathbf{B}\myvec{u},
\end{align}
and then, by taking $\myvec{u}:=\mathbf{B}_2^{-1}\myvec{v}$,
\begin{align}\label{lemp: Theta formula}
\myvec{v}^T\kbra{\mathbf{\Theta}-(\mathbf{B}_2^{-1})^T\mathbf{B}\mathbf{B}_2^{-1}}\myvec{v}\equiv0\quad
\text{or}\quad \mathbf{\Theta}=(\mathbf{B}_2^{-1})^T\mathbf{B}\mathbf{B}_2^{-1}.
\end{align}
By virtue of the decomposition
$\mathbf{B}=\bra{\mathbf{L}\Lambda_\tau^{-1}}^T\mathbf{L}\Lambda_\tau^{-1}$ in Lemma \ref{lem: tilde B-positiveDefinite},
we have
\begin{align*}
\sum_{k=1}^n  \sum_{j=1}^k\theta_{k-j}^{(k)} w_k v_j
=&\,\sum_{k=1}^nw_ku_k=\myvec{u}^T\myvec{w}
=\brab{\mathbf{L}\Lambda_\tau^{-1}\myvec{u}}^T
\brab{\Lambda_\tau\mathbf{L}^{-1}}^T\myvec{w}\nonumber\\
\le&\,\frac{\varepsilon}{2}\myvec{u}^T\brab{\mathbf{L}\Lambda_\tau^{-1}}^T\mathbf{L}\Lambda_\tau^{-1}\myvec{u}
+\frac{1}{2\varepsilon}\myvec{w}^T\brab{\Lambda_\tau\mathbf{L}^{-1}}\brab{\Lambda_\tau\mathbf{L}^{-1}}^T\myvec{w}\nonumber\\
=&\,\frac{\varepsilon}{2}\myvec{u}^T\mathbf{B}\myvec{u}
+\frac{1}{2\varepsilon}\myvec{w}^T\Lambda_\tau\mathbf{L}^{-1}\brab{\mathbf{L}^{-1}}^T\Lambda_\tau\myvec{w}\nonumber\\
=&\,\frac{\varepsilon}{2}\myvec{v}^T\mathbf{\Theta}\myvec{v}
+\frac{1}{2\varepsilon}\myvec{w}^T\mathbf{B}^{-1}\myvec{w}\quad\text{for any $\varepsilon>0$,}
\end{align*}
where the Young's inequality was used in the inequality and the identity \eqref{lemp:conv identity} was applied in the last equality.
This completes the proof.
\end{proof}

Now we are in position to present the proof of Lemma \ref{lem:quadr form inequ}.

\begin{proof}\textbf{\!\!of Lemma \ref{lem:quadr form inequ}}\quad
To avoid possible confusions, we define the vector norm $\timenorm{\cdot}$ by
$\timenorm{\myvec{u}}:=\sqrt{\myvec{u}^T\myvec{u}}$ and
the associated matrix norm
$\timenorm{\textbf{A}}:=\sqrt{\rho\brab{\textbf{A}^T\textbf{A}}}$.
Lemma \ref{lem:quadr form inequ B} gives
\begin{align}\label{lemp:conv Young ineq}
\sum_{k=1}^n  \sum_{j=1}^k\theta_{k-j}^{(k)} w_k v_j
\le
\varepsilon\sum_{k=1}^n  \sum_{j=1}^k\theta_{k-j}^{(k)} v_k v_j
+\frac{1}{2\varepsilon}\myvec{w}^T\mathbf{B}^{-1}\myvec{w}
\quad \text{for $\varepsilon > 0$}.
\end{align}
We will handle the second term
at the right side of \eqref{lemp:conv Young ineq}.
Lemma \ref{lem: tilde B-positiveDefinite} shows $\widetilde{\mathbf{B}}=\mathbf{L}^T\mathbf{L}$ and
$\mathbf{B}^{-1}=\Lambda_\tau\mathbf{L}^{-1}\brab{\mathbf{L}^{-1}}^T\Lambda_\tau$. Moreover,
Lemma \ref{lem:quadr form inequ B} gives
 \begin{align*}
\mathbf{\Theta}=&\,(\mathbf{B}_2^{-1})^T\mathbf{B}\mathbf{B}_2^{-1}
 =(\mathbf{B}_2^{-1})^T\bra{\mathbf{L}\Lambda_\tau^{-1}}^T\mathbf{L}\Lambda_\tau^{-1}\mathbf{B}_2^{-1}
 =\bra{\mathbf{L}\Lambda_\tau^{-1}\mathbf{B}_2^{-1}}^T\mathbf{L}\Lambda_\tau^{-1}\mathbf{B}_2^{-1}
\end{align*}
such that $\myvec{w}^T\mathbf{\Theta}\myvec{w}=\timenorm{\mathbf{L}\Lambda_\tau^{-1}\mathbf{B}_2^{-1}\myvec{w}}^2$.
We apply the definition \eqref{matrix: modified-B2} to derive that
\begin{align*}
\myvec{w}^T\mathbf{B}^{-1}\myvec{w}
=&\,\bra{\brab{\mathbf{L}^{-1}}^T\Lambda_\tau\myvec{w}}^T\brab{\mathbf{L}^{-1}}^T\Lambda_\tau\myvec{w}
=\timenorm{\brab{\mathbf{L}^{-1}}^T\Lambda_\tau\myvec{w}}^2\nonumber\\
=&\,\timenorm{\brab{\mathbf{L}^{-1}}^T\Lambda_\tau\mathbf{B}_2\Lambda_\tau
\mathbf{L}^{-1}\mathbf{L}\Lambda_\tau^{-1}\mathbf{B}_2^{-1}\myvec{w}}^2\nonumber\\
\le&\,\timenorm{\brab{\mathbf{L}^{-1}}^T\Lambda_\tau\mathbf{B}_2
\Lambda_\tau\mathbf{L}^{-1}}^2
\timenorm{\mathbf{L}\Lambda_\tau^{-1}\mathbf{B}_2^{-1}\myvec{w}}^2\\
=&\,\timenorm{\brab{\mathbf{L}^{-1}}^T\widetilde{\mathbf{B}}_2\mathbf{L}^{-1}}^2\cdot
\myvec{w}^T\mathbf{\Theta}\myvec{w}
\le\mathcal{M}_r^{(n)}\cdot
\myvec{w}^T\mathbf{\Theta}\myvec{w},
\end{align*}
where we denote
\begin{align}\label{lemp: quantity M_rn}
\mathcal{M}_r^{(n)}:=&\,\timenorm{\mathbf{L}^{-1}}^4
\timenorm{\widetilde{\mathbf{B}}_2}^2
=\lambda_{\max}^2\brab{\brab{\widetilde{\mathbf{B}}^T}^{-1}}
\lambda_{\max}\brab{\widetilde{\mathbf{B}}_2^T\widetilde{\mathbf{B}}_2}
=\frac{\lambda_{\max}\brab{\widetilde{\mathbf{B}}_2^T\widetilde{\mathbf{B}}_2}}{\lambda_{\min}^{2}\brab{\widetilde{\mathbf{B}}}}.
\end{align}
Therefore it follows from \eqref{lemp:conv Young ineq} that
\begin{align*}
\sum_{k=1}^n  \sum_{j=1}^k\theta_{k-j}^{(k)} w_k v_j
\le
\varepsilon\sum_{k=1}^n  \sum_{j=1}^k\theta_{k-j}^{(k)} v_k v_j
+\frac{\mathcal{M}_r^{(n)}}{\varepsilon}\sum_{k=1}^n  \sum_{j=1}^k\theta_{k-j}^{(k)} w_k w_j
\quad \text{for $\varepsilon > 0$}.
\end{align*}
To complete the proof, it remains to show that $\mathcal{M}_r^{(n)}$ is uniformly bounded with respect to the level index $n$.
Fortunately, Lemmas \ref{lem: tilde B-positiveDefinite} and \ref{lem: tilde B2TB2-positiveDefinite} confirm that
there exists an $n$-independent constant $\mathcal{M}_r:=\max_{n\ge1}\mathcal{M}_r^{(n)}<39$.
Actually, under the weak step-ratio condition \textbf{S1}, one has a rough estimate
\begin{align*}
\mathcal{M}_r=\max_{n\ge1}\frac{\lambda_{\max}\brab{\widetilde{\mathbf{B}}_2^T\widetilde{\mathbf{B}}_2}}{\lambda_{\min}^{2}\brab{\widetilde{\mathbf{B}}}}
\le \max_{n\ge1}\frac{\max_{1\le k\le n}R_U\bra{r_{k},r_{k+1}}}{\min_{1\le k\le n}R_L^2\bra{r_{k},r_{k+1}}}<39.
\end{align*}
It completes the proof.
\end{proof}

\begin{remark}[Improved estimate on the constant $\mathcal{M}_r$]\label{remark: quantity M_rn}
As noted in \cite{Liao2019Adaptive}, the adjacent step ratios take $r_n\approx1$ when the solution varies slowly,
and the restriction \textbf{S1} only takes its effect inside the fast-varying (high gradient) time domains ($r_n<1$),
in the transition regions from the slow-varying to fast-varying domains ($r_n<1$),
and in the ``fast-to-slow'' transition regions ($r_n>1$).
Then the positive constant $\mathcal{M}_r$ in the proof of Lemma \ref{lem:quadr form inequ} can be refined by
considering three different cases,  cf. Remark \ref{remark:comments on time-step condition},
\begin{itemize}
  \item[(a)] If $0<r_{n}\le \sqrt{3}-1$, one can choose
   $\mathcal{M}_r=R_U(\sqrt{3}-1,\sqrt{3}-1)/R_L^{2}\brab{0,\sqrt{3}-1}<1.19;$
  \item[(b)] If $\sqrt{3}-1<r_n\le 2$, then one has
  $\mathcal{M}_r=R_U(2,2)/R_L^{2}\brab{2,2}<3.25;$
  \item[(c)] If $2<r_{n}< r_{\mathrm{sup}}$,  one can choose the next step-ratio $0<r_{n+1}\le 1.45$
  (to reduce or slightly enlarge the step size)
   such that $\mathcal{M}_r= R_U(r_{\mathrm{sup}},1.45)/R_L^{2}\bra{r_{\mathrm{sup}},1.45}<3.94.$
\end{itemize}
Always, one can take $\mathcal{M}_r=4$ in the adaptive computations.
\end{remark}


\end{document}